\documentclass{amsart}
\usepackage{array}

\def\C{\mathbb C}
\def\P{\mathbb P}
\def\R{\mathbb R}
\def\N{\mathbb N}
\def\Z{\mathbb Z}
\def\Q{\mathbb Q}
\def\I{\mathbf I}
\def\Sp{\mathrm{Sp}}
\def\Re{\mathrm{Re\,}}
\newtheorem{Theorem}{Theorem}
\newtheorem{Lemma}{Lemma}

\newtheorem{Remark}{Remark}
\def\Im{\mathrm{Im\,}}

\def\MF#1#2#3#4#5{\begin{pmatrix} \MMF #1 \\ \MMF #2 \\ \MMF #3 \\
  \MMF #4 \\ \MMF #5\end{pmatrix}}
\def\MMF#1#2#3#4#5{#1 & #2 & #3 & #4 & #5}

\def\M#1#2#3#4{\text{\tiny $\begin{pmatrix} \MM #1 \\ \MM #2 \\ \MM #3 \\
  \MM #4 \end{pmatrix}$}}
\def\MM#1#2#3#4{#1 & #2 & #3 & #4}

\begin{document}
\title[Monodromy of Calabi-Yau differential equations]{Monodromy of
  Picard-Fuchs differential equations for Calabi-Yau threefolds}
\author{Yao-Han Chen}
\address{Department of Applied Mathematics \\
  National Chiao Tung University \\
  Hsinchu 300 \\
  TAIWAN}
\email{peace.am92g@nctu.edu.tw}
\author{Yifan Yang}
\address{Department of Applied Mathematics \\
  National Chiao Tung University \\
  Hsinchu 300 \\
  TAIWAN}
\email{yfyang@math.nctu.edu.tw}
\author{Noriko Yui}
\address{Department of Mathematics and Statistics\\
  Queen's University\\
  Kingston, Ontario K3L 3N6\\
  CANADA}
\email{yui@mast.queensu.ca}
\date{16 January 2007}
\thanks{Y.-H. Chen and Yifan Yang were supported by Grant
  94-2115-M-009-012 of the National Science Council (NSC) of the
  Republic of China (Taiwan). N. Yui was supported in part by
  Discovery Grant of the Natural Sciences and Engineering Research 
  Council (NSERC) of Canada}



\begin{abstract} In this paper we are concerned with the monodromy of
  Picard-Fuchs differential equations associated with one-parameter
  families of Calabi-Yau threefolds. Our results show that in the
  hypergeometric cases the matrix representations of monodromy relative
  to the Frobenius bases can be expressed in terms of the geometric
  invariants of the underlying Calabi-Yau threefolds. This phenomenon
  is also verified numerically for other families of Calabi-Yau
  threefolds in the paper. Furthermore, we discover that under a
  suitable change of bases the monodromy groups are contained in
  certain congruence subgroups of $\Sp(4,\Z)$ of finite index and 
  whose levels are related to the geometric invariants of the
  Calabi-Yau threefolds.
\end{abstract}
\maketitle

\begin{section}{Introduction} Let $M_z$ be a family of Calabi-Yau
  $n$-folds parameterized by a complex variable $z\in\P^1(\C)$, and
  $\omega_z$ be the unique holomorphic differential $n$-form on $M_z$
  (up to a scalar). Then the standard theory of Gauss-Manin
  connections asserts that the periods
  $$
    \int_{\gamma_z}\omega_z
  $$
  satisfy certain linear differential equations, called the {\it
  Picard-Fuchs differential equations}, where $\gamma_z$ are
  $r$-cycles on $M_z$.

  When $n=1$, Calabi-Yau onefolds are just elliptic curves. A
  classical example of Picard-Fuchs differential equations is
  \begin{equation} \label{elliptic example}
    (1-z)\theta^2f-z\theta f-\frac z4f=0,\qquad \theta=zd/dz,
  \end{equation}
  satisfied by the periods
  $$
    f(z)=\int^\infty_1\frac{dx}{\sqrt{x(x-1)(x-z)}}
  $$
  of the family of elliptic curves $E_z:y^2=x(x-1)(x-z)$.

  When $n=2$, Calabi-Yau manifolds are either $2$-dimensional complex
  tori or $K3$ surfaces. When the Picard number of a one-parameter
  family of $K3$ surfaces is $19$, the Picard-Fuchs differential
  equation has order $3$. One of the simplest examples is
  $$
    x_1^4+x_2^4+x_3^4+x_4^4-z^{-1}x_1x_2x_3x_4=0\subset\P^3,
  $$
  whose Picard-Fuchs differential operator is
  \begin{equation} \label{K3 example}
    \theta^3-4z(4\theta+1)(4\theta+2)(4\theta+3).
  \end{equation}
  Another well-known example is
  \begin{equation} \label{Apery example}
    (1-34z+z^2)\theta^3+(3z^2-51z)\theta^2+(3z^2-27z)\theta+(z^2-5z),
  \end{equation}
  which is the Picard-Fuchs differential operator for the family of $K3$
  surfaces
  $$
    1-(1-XY)Z-zXYZ(1-X)(1-Y)(1-Z)=0.
  $$
  (See \cite{Beukers-Peters}.) This differential equation appeared in
  Ap\'ery's proof of irrationality of $\zeta(3)$. (See \cite{Beukers}.)

  When $n=3$ and Calabi--Yau threefolds have the Hodge number $h^{2,1}$ 
  equal to $1$, the Picard-Fuchs differential equations have order $4$.
  One of the most well-known examples of such Calabi--Yau threefolds
  is the quintic threefold
  $$
    x_1^5+x_2^5+x_3^5+x_4^5+x_5^5-z^{-1}x_1x_2x_3x_4x_5=0\subset\P^4.
  $$
  In \cite{Candelas}, it is shown that the Picard-Fuchs
  differential operator for this family of Calabi--Yau threefolds is
  \begin{equation} \label{quintic example}
    \theta^4-5z(5\theta+1)(5\theta+2)(5\theta+3)(5\theta+4).
  \end{equation}
  (Actually, it is the mirror partner of
  the quintic Calabi--Yau threefolds that has Hodge number $h^{2,1}=1$
  and hence the Picard--Fuch differential equation is of order
  $4$. But the mirror pair of Calabi--Yau threefolds share the same
  ``principle periods''. This means that the Picard--Fuchs
  differential equation of the original quintic Calabi--Yau threefold
  of order $204$ contains the above order $4$ equation as a factor and
  the factors corresponding to the remaining $200$ ``semiperiods''.

  In this article we are concerned with the monodromy aspect of the
  Picard-Fuchs differential equations. Let
  $$
    L:\,r_n(z)\theta^n+r_{n-1}(z)\theta^{n-1}+\cdots+r_0(z), \qquad
    r_i\in\C(z),
  $$
  be a differential operator with regular singularities. Let $z_0$ be
  a singular point and $S$ be the solution space of $L$ at $z_0$. Then
  analytic continuation along a closed curve $\gamma$ circling $z_0$
  gives rise to an automorphism of $S$, called {\it monodromy}. If a
  basis $\{f_1,\ldots,f_n\}$ of $S$ is chosen, then we have a matrix
  representation of the monodromy. Suppose that $f_i$ becomes
  $a_{i1}f_1+\cdots+a_{in}f_n$ after completing the loop $\gamma$,
  that is, if
  $$
    \begin{pmatrix} f_1 \\ \vdots \\ f_n\end{pmatrix} \longmapsto
    \begin{pmatrix} a_{11} & \ldots & a_{1n} \\
      \vdots &  & \vdots \\ a_{n1} & \ldots & a_{nn} \end{pmatrix}
    \begin{pmatrix} f_1 \\ \vdots \\ f_n\end{pmatrix},
  $$
  then the matrix representation of the monodromy along $\gamma$
  relative to the basis $\{f_i\}$ is the matrix $(a_{ij})$. The group
  of all such matrices are referred to as the {\it monodromy group}
  relative to the basis $\{f_i\}$ of the differential equation.
  Clearly, two different choices of bases may result in two different
  matrix representations for the same monodromy. However, it is easily
  seen that they are connected by conjugation by the matrix of basis
  change. Thus, the monodromy group is defined up to conjugation.
  In the subsequent discussions, for the ease of exposition, 
  we may often drop the phrase ``up to
  conjugation'' about the monodromy groups, 
  when there is no danger of ambiguities.

  It is known that for one-parameter families of Calabi-Yau varieties
  of dimension one and two (i.e., elliptic curves and $K3$ surfaces, 
  respectively), the monodromy groups are very often congruence
  subgroups of $SL(2,\R)$. For instance, the monodromy group of
  (\ref{elliptic example}) is $\Gamma(2)$, while those of (\ref{K3
  example}) and (\ref{Apery example}) are $\Gamma_0(2)+\omega_2$ and
  $\Gamma_0(6)+\omega_6$, respectively, where $\omega_d$ denotes the
  Atkin-Lehner involution. (Technically speaking, the monodromy groups
  of (\ref{K3 example}) and (\ref{Apery example}) are subgroups of
  $SL(3,\R)$ since the order of the differential equations is $3$.
  But because (\ref{K3 example}) and (\ref{Apery example}) are
  symmetric squares of second-order differential equations, we may
  describe the monodromy in terms of the second-order ones.)
  Moreover, suppose that $y_0(z)=1+\cdots$ is
  the unique holomorphic solution at $z=0$ and $y_1(z)=y_0(z)\log
  z+g(z)$ is the solution with logarithmic singularity. Set
  $\tau=cy_1(z)/y_0(z)$ for a suitable complex number $c$. Then $z$,
  as a function of $\tau$, becomes a modular function, and
  $y_0(z(\tau))$ becomes a modular form of weight $1$ for the order
  $2$ cases and of weight $2$ for the order $3$ cases. For example, a
  classical result going back to Jacobi states that
  $$
    \theta_3^2=\, _2F_1\left(\frac12,\frac12;1;
    \frac{\theta_2^4}{\theta_3^4}\right),
  $$
  where
  $$
    \theta_2(\tau)=q^{1/8}\sum_{n\in\Z}q^{n(n+1)/2}, \qquad
    \theta_3(\tau)=\sum_{n\in\Z}q^{n^2/2}, \qquad q=e^{2\pi i\tau},
  $$
  or equivalently, that the modular form $y(\tau)=\theta_3^2$, as a
  function of $z(\tau)=\theta_2^4/\theta_3^4$, satisfies
  (\ref{elliptic example}). Here $_2F_1$ denotes the Gauss
  hypergeometric function.

  In this paper we will address the monodromy problem for Calabi-Yau
  threefolds. At first, given the experience with the elliptic curve
  and $K3$ surface cases, one may be tempted to guess that the
  monodromy group of such a differential equation will be the
  symmetric cube of some congruence subgroup of $SL(2,\R)$. After all,
  there is a result by Stiller \cite{Stiller} (see also \cite{Yang})
  asserting that if $t(\tau)$ is a non-constant modular function and
  $F(\tau)$ is a modular form of weight $k$ on a subgroup of
  $SL(2,\R)$ commensurable with $SL(2,\Z)$, then $F,\tau
  F,\ldots,\tau^k F$, as functions of $t$, are solutions of a
  $(k+1)$-st order linear differential equation with algebraic
  functions of $t$ as coefficients. However, this is not the case in
  general. A quick way to see this is that the coefficients of the
  symmetric cube of a second order differential equation
  $y^{\prime\prime}+r_1(t)y^\prime+r_0(t)y=0$ is completely determined
  by $r_1$ and $r_0$, but the coefficients of the Picard-Fuchs
  differential equations, including (\ref{quintic example}), do not
  satisfy the required relations. (The exact relations can be computed
  using Maple's command {\tt symmetric\_power}.) Nevertheless, in the
  subsequent discussion we will show that, with a suitable choice of
  bases, the monodromy groups for Calabi-Yau threefolds are contained
  in certain congruence subgroups of $\Sp(4,\Z)$ whose levels are
  somehow described in terms of the geometric invariants of the 
  manifolds in question. This is proved rigorously for the
  hypergeometric cases and verified numerically for other (e.g.,
  non-hypergeometric) cases. Furthermore, our
  computation in the hypergeometric cases shows that the matrix
  representation of the monodromy around the finite singular point 
  (different from the origin) relative to the Frobenius basis at 
  the origin can be expressed completely using the geometric invariants of the
  associated Calabi-Yau threefolds. This phenomenon is also verified
  numerically in the non-hypergeometric cases. 
  Although it is highly expected that
  geometric invariants will enter into the picture, in reality,
  geometry will dominate the entire picture in the sense that 
  {\it every} entry of the matrix is expressed
  exclusively in terms of the geometric invariants.

  The monodromy problem in general has been addressed by a number of 
  authors. Papers relevant to our consideration include
  \cite{Beukers-Heckman}, \cite{Candelas}, \cite{Doran-Morgan},
  \cite{Klemm-Theisen1}, and \cite{Enckevort-Straten}, to name a few.
  In \cite{Beukers-Heckman}, Beukers and Heckman studied monodromy
  groups for the hypergeometric functions $_n F_{n-1}$. They 
  showed that the Zariski closure of the monodromy groups of
  (\ref{quintic example}) is $\Sp(4,\C)$. The same is true for other
  Picard-Fuchs differential equations for Calabi-Yau threefolds that
  are hypergeometric. In \cite{Candelas}, Candelas et al. obtained
  precise matrix representations of monodromy for (\ref{quintic
  example}). Then Klemm and Theisen \cite{Klemm-Theisen1} applied the
  same method as that of Candelas et al. to deduce monodromy groups 
  for three other hypergeometric cases. In \cite{Doran-Morgan} Doran 
  and Morgan determined the monodromy groups for all the hypergeometric cases.
  Their matrix representations also involve geometric invariants of
  the Calabi-Yau threefolds. For Picard-Fuchs differential equations
  of non-hypergeometric type, there is not much known in literature.
  In \cite{Enckevort-Straten} van Enckevort and van Straten computed
  the monodromy matrices numerically for a large class of differential
  equations. In many cases, they are able to find bases such that the
  monodromy matrices have rational entries. We will discuss the above
  results in more detail in Sections 3--5.

  Our motivations of this paper may be formulated as follows.   
  Modular functions and modular forms have been extensively investigated 
  over the years, and there are great body of literatures on these subjects.
  As we illustrated above, the monodromy groups of Picard-Fuchs
  differential equations for families of elliptic curves and $K3$
  surfaces are congruence subgroups of $SL(2,\R)$. This modularity
  property can be used to
  study properties of the differential equations and the associated
  manifolds. For instance, in \cite{Lian-Yau1} Lian and Yau gave a uniform
  proof of the integrality of Fourier coefficients of the mirror maps
  for several families of $K3$ surfaces using the fact that the
  monodromy groups are congruence subgroups of $SL_2(\R)$. For such an
  application, it is important to express monodromy groups in a proper
  way so that properties of the associated differential equations can
  be more easily discussed and obtained. Thus, the main motivation of
  our investigation is to find a good representation for monodromy
  groups from which further properties of Picard-Fuchs
  differential equations for Calabi-Yau threefolds can be derived.
  
  The terminology ``modularity'' has been used for many
  different things. One aspect of the modularity that we would like to 
  address is the modularity question of the Galois representations attached 
  to Calabi--Yau threefolds, assuming that Calabi--Yau threefolds in
  question are defined over $\Q$. Let $X$ be a Calabi--Yau threefold
  defined over $\Q$. We consider the $L$-series associated to the 
  third \'etale cohomology group of $X$. It is expected that the
  $L$-series should be determined by some modular (automorphic) forms.
  The examples of Calabi--Yau threefolds we treat in this paper are those
  with the third Betti number equal to $4$.
  It appears that Calabi--Yau threefolds with this property are
  rather scarce.  Batyrev and Straten \cite{Batyrev-Straten} considered
  $13$ examples of Calabi--Yau threefolds with Picard number
  $h^{1,1}=1$. Then their mirror partners will fulfill this
  requirement. (We note that more examples of such Calabi--Yau
  threefolds were found by Borcea \cite{Borcea}.) All these $13$
  Calabi--Yau threefolds are defined as complete intersections of
  hypersurfaces in weighted projective spaces, and they have defining
  equations defined over $\Q$.

  To address the modularity, we ought to have some ``modular
  groups'', and this paper offers candidates for appropriate
  modular groups via
  the monodromy group of the associated Picard--Fuchs differential
  equation (of order $4$). In these cases, we expect that modular forms
  of more variables, e.g., Siegel modular forms associated to the
  modular groups for our congruence subgroups would enter the scene.

  In general, the third Betti numbers of Calabi--Yau threefolds are
  rather large, and consequently, the dimension of the 
  associated Galois representations would be rather high. To remedy
  this situation, we first decompose Calabi--Yau threefolds into
  motives, and then consider the motivic Galois representations
  and their modularity. Especially, when the principal 
  motives (e.g., the motives that are invariant under the mirror maps)  
  are of dimension $4$, the modularity question for such motives should be  
  accessible using the method developed for the examples discussed in this
  paper.

  The modularity questions will be treated in subsequent papers. 
\end{section}

\begin{section}{Statements of results}
  To state our first result, let us recall that among all the
  Picard-Fuchs differential equations for Calabi-Yau threefolds, there
  are $14$ equations that are hypergeometric of the form
  $$
    \theta^4-Cz(\theta+A)(\theta+1-A)(\theta+B)(\theta+1-B).
  $$
  Their geometric descriptions and references are given in the
  following Table 1.
  $$
    \extrarowheight3pt
    \arraycolsep4pt
    \begin{array}{|c||c|c|c||l|c|c|c||c|} \hline
    \# & A & B & C & \text{Description} & H^3 & c_2\!\cdot\! H\!\! & c_3
      & \text{Ref} \\ \hline
    1 & 1/5 & 2/5 & 3125 & X(5)\subset\P^4
      & 5 & 50 & -200 & \text{\cite{Candelas}} \\ \hline
    2 & 1/10 & 3/10 & 8\cdot 10^5 & X(10)\subset\P^4(1,1,1,2,5)
      & 1 & 34 & -288 & \text{\cite{Morrison}} \\ \hline
    3 & 1/2 & 1/2 & 256 & X(2,2,2,2)\subset\P^7
      & 16 & 64 & -128 & \text{\cite{Libgober-Teitelbaum}} \\ \hline
    4 & 1/3 & 1/3 & 729 & X(3,3)\subset\P^5
      & 9 & 54 & -144 & \text{\cite{Libgober-Teitelbaum}} \\ \hline
    5 & 1/3 & 1/2 & 432 & X(2,2,3)\subset\P^6
      & 12 & 60 & -144 & \text{\cite{Libgober-Teitelbaum}} \\ \hline
    6 & 1/4 & 1/2 & 1024 & X(2,4)\subset\P^5
      & 8 & 56 & -176 & \text{\cite{Libgober-Teitelbaum}} \\ \hline
    7 & 1/8 & 3/8 & 65536 & X(8)\subset\P^4(1,1,1,1,4)
      & 2 & 44 & -296 & \text{\cite{Morrison}} \\ \hline
    8 & 1/6 & 1/3 & 11664 & X(6)\subset\P^4(1,1,1,1,2)
      & 3 & 42 & -204 & \text{\cite{Morrison}} \\ \hline
    9 & 1/12 & 5/12 & 12^6 & X(2,12)\subset\P^5(1,1,1,1,4,6)
      & 1 & 46 & -484 & \text{\cite{Doran-Morgan}} \\ \hline
   10 & 1/4 & 1/4 & 4096 & X(4,4)\subset\P^5(1,1,1,1,2,2)
      & 4 & 40 & -144 & \text{\cite{Klemm-Theisen2}} \\ \hline
   11 & 1/4 & 1/3 & 1728 & X(4,6)\subset\P^5(1,1,1,2,2,3)
      & 6 & 48 & -156 & \text{\cite{Klemm-Theisen2}} \\ \hline
   12 & 1/6 & 1/4 & 27648 & X(3,4)\subset\P^5(1,1,1,1,1,2)
      & 2 & 32 & -156 & \text{\cite{Klemm-Theisen2}} \\ \hline
   13 & 1/6 & 1/6 & 2^8\cdot 3^6 & X(6,6)\subset\P^5(1,1,2,2,3,3)
      & 1 & 22 & -120 & \text{\cite{Klemm-Theisen2}} \\ \hline
   14 & 1/6 & 1/2 & 6912 & X(2,6)\subset\P^5(1,1,1,1,1,3)
      & 4 & 52 & -256 & \text{\cite{Klemm-Theisen2}} \\ \hline
    \end{array}
  $$  
  Some comments might be in order for the notations in the table. We
  employ the notations of van Enckevort and van Straten
  \cite{Enckevort-Straten}.
  $X(d_1,d_2,\ldots,d_k)\subset\P^n(w_0,\ldots,w_n)$ stands for a
  complete intersection of $k$ hypersurfaces of degree
  $d_1,\ldots,d_k$ in the weighted projective space with weight
  $(w_0,\cdots, w_n)$. For instance, $X(3,3)\subset\P^5$ is a complete
  intersection of two cubics in the ordinary projective $5$-space
  $\P^5$ defined by
  $$\{\,Y_1^3+Y_2^3+Y_3^3-3\phi Y_4Y_5Y_6=0\,\}\cap\{\,
                 -3\phi Y_1Y_2Y_3+Y_4^3+Y_5^3+Y_6^3=0\,\}.
  $$
  Slightly more generally, $X(4,4)\subset\P^5(1,1,2,1,1,2)$ denotes a
  complete intersection of two quartics in the weighted projective
  $5$-space $\P^5(1,1,2,1,1,2)$ and may be defined by the equations
  $$\{\,Y_1^4+Y_2^4+Y_3^2-4\phi Y_4Y_5Y_6=0\,\}\cap\{\,
                 -4\phi Y_1Y_2Y_3 + Y_4^4+Y_5^4+Y_6^2=0\,\}.
  $$
  We note that all these examples of Calabi--Yau threefolds $M$ have
  the Picard number $h^{1,1}(M)=1$.  Let $\mathcal O(H)$ be the ample
  generator of the Picard group $Pic(M)\simeq \Z$.
  The basic invariants for such a Calabi--Yau threefold $M$ are
  the {\it degree} $d:=H^3$, the {\it second Chern number} $c_2\cdot
  H$ and the Euler number $c_3$ (the Euler characteristic of $M$). The
  equations are numbered in the same way as in \cite{AESZ}.

  In \cite{Candelas}, using analytic properties of hypergeometric
  functions, Candelas et al. proved that with 
  respect to a certain basis, the monodromy matrices around $z=0$ and
  $z=1/3125$ for the quintic threefold case (Equation 1 from Table 1) are
  $$
    \begin{pmatrix} 51 & 90 & -25 & 0 \\
      0 & 1 & 0 & 0 \\ 100 & 175 & -49 & 0 \\
      -75 & -125 & 35 & 1 \end{pmatrix} \qquad \text{and}\qquad
    \begin{pmatrix} 1 & 0 & 0 & 0 \\ 0 & 1 & 0 & 1 \\
    0 & 0 & 1 & 0 \\ 0 & 0 & 0 & 1 \end{pmatrix},
  $$
  respectively. (Note that these two matrices are both in $\Sp(4,\Z)$.)
  Applying the same method as that of Candelas et al., 
  Klemm and Theisen \cite{Klemm-Theisen1}
  also obtained the monodromy of the one-parameter families of
  Calabi--Yau threefolds for Equations 2, 7, and 8. Presumably,
  their method should also work for several other hypergeometric
  cases. However, the method fails when the indicial equation of the
  singularity $\infty$ has repeated roots. To be more precise, it does
  not work for Equations 3--6, 10, 13 and 14.
  Moreover, the method uses the explicit knowledge that the singular
  point $z=1/C$ is of conifold type. (Note that in geometric terms, a
  conical singularity is a regular singular point whose neighborhood
  looks like a cone with a certain base. For instance, a
  $3$-dimensional conifold singularity is locally isomorphic to
  $XY-ZT=0$ or equivalently, to $X^2+Y^2+Z^2+T^2=0$. Reflecting to the
  Picard-Fuchs differential equations, this means that the local
  monodromy is unipotent of index $1$.) Thus, it can not be applied
  immediately to study monodromy of general hypergeometric
  differential equations.

  In \cite{Doran-Morgan} Doran and Morgan proved that if the
  characteristic polynomial of the monodromy around $\infty$ is
  $$
    x^4+(k-4)x^3+(6-2k+d)x^2+(k-4)x+1,
  $$
  then there is a basis such that the monodromy matrices around $z=0$
  and $z=1/C$ are
  \begin{equation} \label{DM representation}
    \begin{pmatrix} 1 & 1 & 0 & 0 \\ 0 & 1 & d & 0 \\
    0 & 0 & 1 & 1 \\ 0 & 0 & 0 & 1\end{pmatrix} \qquad\text{and}\qquad
    \begin{pmatrix} 1 & 0 & 0 & 0 \\ -k & 1 & 0 & 0 \\
    -1 & 0 & 1 & 0 \\ -1 & 0 & 0 & 1 \end{pmatrix},
  \end{equation}
  respectively. It turns out that these numbers $d$ and $k$ both
  have geometric interpretation. Namely, the number $d=H^3$ is the
  degree of the associated threefolds and $k=c_2\!\cdot\! H/12+H^3/6$ is
  the dimension of the linear system $|H|$. Doran and Morgan's
  representation has the advantage that the geometric invariants can
  be read off from the matrices directly (although there is no way to
  extract the Euler number $c_3$ from the matrices), but has the
  disadvantage that the matrices are no longer in the symplectic group
  (in the strict sense).

  Before we state our Theorem 1, let us recall the definition of
  Frobenius basis. Since the only solution of the indicial equation at
  $z=0$ for each of the cases is $0$ with multiplicity $4$, the
  monodromy around $z=0$ is maximally unipotent. (See \cite{Morrison}
  for more detail.) Then the standard method of Frobenius implies
  that at $z=0$ there are four solutions $y_j$, $j=0,\ldots,$, with
  the property that
  \begin{equation} \label{Frobenius basis}
  \begin{split}
    y_0=1+\cdots, & \quad y_1=y_0\log z+g_1, \\
    y_2=\frac12 y_0\log^2 z+g_1\log z+g_2, &\quad
    y_3=\frac16 y_0\log^3 z+\frac12g_1\log^2 z+g_2\log z+g_3,
  \end{split}
  \end{equation}
  where $g_i$ are all functions holomorphic and vanishing at $z=0$.
  We remark that these solutions satisfy the relation
  $$
    \begin{vmatrix} y_0 & y_3 \\ y_0^\prime & y_3^\prime\end{vmatrix}
   =\begin{vmatrix} y_1 & y_2 \\ y_1^\prime & y_2^\prime\end{vmatrix},
  $$
  and therefore the monodromy matrices relative to the ordered basis
  $\{y_0,y_2,y_3,y_1\}$ are in $\Sp(4,\C)$, as predicted by
  \cite{Beukers-Heckman}. Now we can present our first theorem. 

\begin{Theorem} Let
  $$
    L: \theta^4-Cz(\theta+A)(\theta+1-A)(\theta+B)(\theta+1-B)
  $$
  be one of the $14$ hypergeometric equations, and $H^3$, $c_2\!\cdot\!
  H$, and $c_3$ be geometric invariants of the associated Calabi-Yau
  threefolds given in the table above. Let $y_j$, $j=0,\ldots,3$, be
  the Frobenius basis specified by (\ref{Frobenius basis}). Then with
  respect to the ordered basis $\{y_3/(2\pi i)^3,y_2/(2\pi i)^2,
  y_1/(2\pi i), y_0\}$, the monodromy matrices around $z=0$ and
  $z=1/C$ are
  \begin{equation} \label{raw matrices}
    \begin{pmatrix} 1 & 1 & 1/2 & 1/6 \\ 0 & 1 & 1 & 1/2 \\
     0 & 0 & 1 & 1 \\ 0 & 0 & 0 & 1 \end{pmatrix}\quad\text{and}\quad
    \begin{pmatrix} 1+a & 0 & ab/d & a^2/d \\
    -b & 1 & -b^2/d & -ab/d \\ 0 & 0 & 1 & 0 \\
    -d & 0 & -b & 1-a \end{pmatrix},
  \end{equation}
  respectively, where
  $$
    a=\frac{c_3}{(2\pi i)^3}\zeta(3), \qquad b=c_2\!\cdot\! H/24,
    \qquad d=H^3.
  $$
\end{Theorem}

\begin{Remark}
We remark that by conjugating by the matrix
$$
  \begin{pmatrix} d & 0 & b & a \\ 0 & d & d/2 & d/6+b \\
  0 & 0 & 1 & 1 \\ 0 & 0 & 0 & 1 \end{pmatrix},
$$
we do recover Doran and Morgan's representation (\ref{DM
  representation}).  Thus, our Theorem 1 strenthens
the results of Doran and Morgan \cite{Doran-Morgan}. 
Although the referee suggested that Theorem 1 might be a reformulation of the results of 
Doran and Morgan. However, we do not believe that is the case. For one
thing, the argument of Doran--Morgan is purely based on Linear Algebra. It might
be possible to derive our Theorem 1 combining the results of
Doran--Morgan and those of Kontsevich; we will not address this
question here, but left to future investigations.
\end{Remark}

The appearance of the geometric invariants $c_2,\, c_3,\, H$ and $d$ 
is not so surprising. In \cite{Candelas}, it was shown that the conifold 
period, defined up to a constant as the holomorphic solution
$f(z)=a_1(z-1/C)+a_2(z-1/C)^2+\cdots$ at $z=1/C$ that appears in the
unique solution $f(z)\log(z-1/C)+g(z)$ with logarithmic singularity at
$z=1/C$, is asymptotically
\begin{equation} \label{conifold period}
  \frac{H^3}{6(2\pi i)^3}\log^3 z+\frac{c_2\!\cdot\! H}{48\pi i}
  \log z+\frac{c_3}{(2\pi i)^3}\zeta(3)+\cdots
\end{equation}
near $z=0$. (See also \cite{Katz-Klemm-Vafa}.) Therefore, it is
expected that the entries of the monodromy matrices should contain
the invariants. However, it is still quite remarkable that the matrix
is determined completely by the invariants alone. We have
numerically verified the phenomenon for other families of Calabi-Yau
threefolds, and also for general differential equations of Calabi-Yau
type. (See \cite{AESZ} for the definition of a differential equation
of Calabi-Yau type. See also Section 5 below.) It appears that if the
differential equation has at least one singularity with exponents
$0,1,1,2$, then there is always a singularity whose monodromy relative
to the Frobenius basis is of the form stated in the theorem. Thus,
this gives a numerical method to identify the possible geometric
origin of a differential equation of Calabi-Yau type.

We emphasize that our proof of Theorem 1 is merely a verification. That is,
we can prove it, but unfortunately it does
not give any geometric insight why the matrices are in this special
form. 

Acutally, the referee has pointed out that such a geometric
  interpretation seems to exist by Kontsevich.  In the framework of
  ``homological mirror symmetry'' of Kontsevich, the first matrix in
  Theorem 1 would be the matrix associated to tensoring by the
  hyperplane line bundle in the bounded derived category of sheaves on
  the Calabi--Yau variety. In general, the matrices in Theorem 1
  describe the cohomology action of certain Fourier--Mukai
  functors. In particular, this explains why the matrices are
  determined by topological invariants of the underlying Calabi--Yau
  manifolds. The paper of van Enckevort and van Straten
  \cite{Enckevort-Straten} addressed monodromy calculations of fourth
  order equations of Calabi--Yau type based on homological mirror
  symmetry. The reader is referred to the article
  \cite{Enckevort-Straten} for full details about geometric
  interpretations of matrices.  We wonder, though, if the Kontsevich's results 
  fully explain why there are no ``non-geometric'' numbers in the second matrices. 
  To be more precise, here is our question. Since the
  second matrix $M$ is unipotent of rank 1, we know that the rows of
  $M-Id$ are all scalar multiples of a fixed row vector. We probably
  can deduce from Kontsevich's result that the fourth row is
  $(-d,0,-b,-a)$, but why the first three rows of $M-Id$ are $-a/d$,
  $-b/a$, and $0$ times this vector (but not other ``non-geometric''
  scalars?

Now conjugating the matrices by
\begin{equation} \label{CY basis}
  \begin{pmatrix} 0 & 0 & 1 & 0 \\ 0 & 0 & 0 & 1 \\
  0 & d & d/2 & -b \\ -d & 0 & -b & -a \end{pmatrix},
\end{equation}
we can bring the matrices into the symplectic group $\Sp(4,\Z)$.
The results are
$$
  \begin{pmatrix} 1 & 1 & 0 & 0 \\ 0 & 1 & 0 & 0 \\
  d & d & 1 & 0 \\ 0 & -k & -1 & 1 \end{pmatrix}\quad\text{and}
  \quad\begin{pmatrix} 1 & 0 & 0 & 0 \\ 0 & 1 & 0 & 1 \\
  0 & 0 & 1 & 0 \\ 0 & 0 & 0 & 1\end{pmatrix}
$$
for $z=0$ and $z=1/C$, respectively, where $k=2b+d/6$. Since the
monodromy group is generated by these two matrices, we see that the
group is contained in the congruence subgroup $\Gamma(d,\gcd(d,k))$,
where the notation $\Gamma(d_1,d_2)$ with $d_2|d_1$ represents
\begin{equation*}
\begin{split}
  \Gamma(d_1,d_2)&=\left\{\gamma\in\Sp(4,\Z):\gamma\equiv
  \begin{pmatrix} 1 & \ast & \ast & \ast \\
   0 & \ast & \ast & \ast \\ 0 & 0 & 1 & 0 \\
   0 & \ast & \ast & \ast \end{pmatrix} \mod d_1\right\} \\
   &\qquad\qquad\bigcap\left\{\gamma\in\Sp(4,\Z):\gamma\equiv
  \begin{pmatrix} 1 & \ast & \ast & \ast \\
   0 & 1 & \ast & \ast \\ 0 & 0 & 1 & 0 \\
   0 & 0 & \ast & 1 \end{pmatrix} \mod d_2\right\}
\end{split}
\end{equation*}
We remark that the entries of the matrices in $\Gamma(d_1,d_2)$
satisfy certain congruence relations inferred from the symplecticity
of the matrices. To be more explicit, let us recall that the
symplectic group is characterized by the property that
$$
  \gamma=\begin{pmatrix} A & B \\ C & D\end{pmatrix} \in \Sp(2n,\C),
$$
where $A$, $B$, $C$, and $D$ are $n\times n$ blocks, if and only if
$$
  \gamma^{-1}=\begin{pmatrix} D^t & -B^t \\ -C^t & A^t \end{pmatrix}.
$$
Thus, for
$$
  \begin{pmatrix} a_{11} & a_{12} & a_{13} & a_{14} \\
  a_{21} & a_{22} & a_{23} & a_{24} \\
  a_{31} & a_{32} & a_{33} & a_{34} \\
  a_{41} & a_{42} & a_{43} & a_{44} \end{pmatrix}
$$
to be in $\Gamma(d_1,d_2)$, the integers $a_{ij}$ should satisfy the
implicit conditions
$$
  a_{22}a_{44}-a_{24}a_{42}\equiv 1,\quad
  a_{23}\equiv a_{14}a_{22}-a_{12}a_{24}, \quad
  a_{43}\equiv a_{14}a_{42}-a_{12}a_{44} \mod d_1,
$$
and
$$
  a_{12}\equiv -a_{43} \mod d_2.
$$
We now summarize our finding in the following theorem.

\begin{Theorem} Let
$$
  \theta^4-Cz(\theta+A)(\theta+1-A)(\theta+B)(\theta+1-B)
$$
be one of the $14$ hypergeometric equations. Let $y_j$, $j=0,\ldots,3$
be the Frobenius basis. Then relative to the ordered basis 
$$
  \frac{y_1}{2\pi i},\quad y_0,\quad \frac{H^3}{2(2\pi i)^2}y_2+
  \frac{H^3}{4\pi i}y_1-\frac{c_2\!\cdot\! H}{24}y_0, \quad
 -\frac{H^3}{6(2\pi i)^3}y_3-\frac{c_2\!\cdot\!H}{48\pi i}y_1
 -\frac{c_3}{(2\pi i)^3}\zeta(3)y_0,
$$
the monodromy matrices around $z=0$ and $z=1/C$ are
\begin{equation} \label{nice matrices}
  \begin{pmatrix} 1 & 1 & 0 & 0 \\ 0 & 1 & 0 & 0 \\
  d & d & 1 & 0 \\ 0 & -k & -1 & 1 \end{pmatrix}\quad\text{and}
  \quad\begin{pmatrix} 1 & 0 & 0 & 0 \\ 0 & 1 & 0 & 1 \\
  0 & 0 & 1 & 0 \\ 0 & 0 & 0 & 1\end{pmatrix}
\end{equation}
with $d=H^3$, $k=H^3/6+c_2\!\cdot\!H/12$, respectively. They are
contained in the congruence subgroups $\Gamma(d_1,d_2)$ for
the $14$ cases in the table below.
$$
  \extrarowheight3pt
  \begin{array}{||c|cc|cc||c|cc|cc||} \hline
  \# & A & B & d_1 & d_2 & \# & A & B & d_1 & d_2   \\ \hline
  1 & 1/5 & 2/5 & 5 & 5 &   8 & 1/6 & 1/3 & 3 & 1   \\ \hline
  2 & 1/10 & 3/10 & 1 & 1 & 9 & 1/12 & 5/12 & 1 & 1 \\ \hline
  3 & 1/2 & 1/2 & 16 & 8 & 10 & 1/4 & 1/4 & 4 & 4   \\ \hline
  4 & 1/3 & 1/3 & 9 & 3 &  11 & 1/4 & 1/3 & 6 & 1   \\ \hline
  5 & 1/3 & 1/2 & 12 & 1 & 12 & 1/6 & 1/4 & 2 & 1   \\ \hline
  6 & 1/4 & 1/2 & 8 & 2 &  13 & 1/6 & 1/6 & 1 & 1   \\ \hline
  7 & 1/8 & 3/8 & 2 & 2 &  14 & 1/6 & 1/2 & 4 & 1   \\ \hline
  \end{array}
$$
\end{Theorem}
\medskip

\noindent{\bf Remark.} We remark that what we show in Theorem 2 is
merely the fact that the monodromy groups are contained in the congruence
subgroups $\Gamma(d_1,d_2)$. Although the congruence subgroups
$\Gamma(d_1,d_2)$ are of finite index in $\Sp(4,\Z)$ (see the appendix
by Cord Erdenberger for the index formula), the monodromy
groups themselves may not be so. 

At this moment, we cannot say anything
definite about these groups, e.g., their finite indexness. In fact, there
are two opposing speculations (one by the authors, and the other
by Zudilin) about these groups. 
We believe, based on a result (Theorem 13.3) of Sullivan \cite{Sullivan}, that it might
be justified to claim that the monodromy group is an arithmetic subgroup of
the congruence subgroup $\Gamma(d_1,d_2)$, and hence is of finite index.
(Andrey Todorov pointed out to us Sullivan's theorem, though we must
confess that we do not fully understand the paper of Sullivan.) 

As opposed to our belief, Zudilin has indicated to us via e-mail that
a heuristic argument suggests that the monodromy groups are too
``thin'' to be of finite index.
\medskip

It would not be of much significance if the hypergeometric equations
are the only cases where the monodromy groups are contained in
congruence subgroups. Our numerical computation suggests that there
are a number of further examples where the monodromy groups continue to be
contained in congruence subgroups of $\Sp(4,\Z)$. However, the general
picture is not as simple as that for the hypergeometric cases.

As mentioned earlier, our numerical data suggest that the Picard-Fuchs
differential equations for Calabi-Yau threefolds known in literature
all have bases relative to which the monodromy matrices around the
origin and some singular points of conifolds take the form (\ref{raw
  matrices}) described in Theorem 1. Thus, with respect to the basis
given in Theorem 2, the matrices around the origin and the conifold
points again have the form (\ref{nice matrices}). However, with this
basis change, the monodromy matrices around other singularities may
not be in $\Sp(4,\Z)$, but in $\Sp(4,\Q)$ instead, although the
entries still satisfy certain congruence relations. Furthermore, in
most cases, we can still realize the
monodromy groups in congruence subgroups of $\Sp(4,\Z)$, by a suitable
conjugation.
\medskip

\noindent{\bf Example 1.} Consider the differential equation
\begin{equation*}
\begin{split}
  &25\theta^4-15z(51\theta^4+84\theta^3+72\theta^2+30\theta+5) \\
  &\qquad  +6z^2(531\theta^4+828\theta^3+541\theta^2+155\theta+15) \\
  &\qquad-54z^3(423\theta^4+2160\theta^3+4399\theta^2+3795\theta+1170)\\
  &\qquad+243z^4(279\theta^4+1368\theta^3+2270\theta^2+1586\theta+402)
  -59049z^5(\theta+1)^4.
\end{split}
\end{equation*}
In \cite{Batyrev-Straten} it is shown that this is the Picard-Fuchs
differential equation for the Calabi-Yau threefolds defined as the
complete intersection of three hypersurfaces of degree $(1,1,1)$ in
$\P^2\times\P^2\times\P^2$. The invariants are $H^3=90$,
$c_2\!\cdot\!H=108$, and $c_3=-90$. There are $6$ singularities $0$,
$1/27$, $\pm i/\sqrt{27}$, $5/9$, and $\infty$ for the differential
equation. Among them, the local exponents at $z=5/9$ are
$0,1,3,4$ and we find that the monodromy around $z=5/9$ is the identity.
For others, our numerical computation shows that relative to the basis
in Theorem 2 the monodromy matrices are
$$
  T_0=\begin{pmatrix}1 & 1 & 0 & 0 \\ 0 & 1 & 0 & 0 \\
  90 & 90 & 1 & 0 \\ 0 & -24 & -1 & 1\end{pmatrix}, \qquad
  T_{1/27}=\begin{pmatrix} 1 & 0 & 0 & 0 \\ 0 & 1 & 0 & 1 \\
   0 & 0 & 1 & 0 \\ 0 & 0 & 0 & 1\end{pmatrix},
$$
$$
  T_{i/\sqrt{27}}=\begin{pmatrix}-17 & 3 & 1/3 & 1\\
  -54 & 10 & 1 & 3 \\ -972 & 162 & 19 & 54 \\
  162 & -27 & -3 & -8 \end{pmatrix}, \qquad
  T_{-i/\sqrt{27}}=\begin{pmatrix} -11 & 3 & 1/3 & -1 \\
   36 & -8 & -1 & 3 \\ -432 & 108 & 13 & -36 \\
  108 & -27 & -3 & 10 \end{pmatrix}.
$$
From these, we see that the monodromy group is contained in the
following group
\begin{equation*}
\begin{split}
 &\Big\{(a_{ij})\in\Sp(4,\Q): a_{ij}\in\Z\ \forall(i,j)\neq(1,3),\
  a_{13}\in\frac13\Z, \\
 &\qquad\ a_{21}, a_{31}, a_{41}, a_{32}, a_{34}\equiv 0\mod 18,\quad
  a_{11},a_{33}\equiv 1\mod 6, \\
 &\qquad\ a_{42}\equiv0, \quad a_{22}, a_{44}\equiv 1\mod 3\Big\}
\end{split}
\end{equation*}
Conjugating by
$$
  \begin{pmatrix} 3 & 0 & 0 & 0 \\ 0 & 1 & 0 & 0 \\
  0 & 0 & 1 & 0 \\ 0 & 0 & 0 & 1 \end{pmatrix},
$$
we find that the monodromy group can be brought into the congruence
subgroup $\Gamma(6,3)$.
\medskip

\noindent{\bf Example 2.} Consider the differential equation.
\begin{equation*}
\begin{split}
  &9\theta^4-3z(173\theta^4+340\theta^3+272\theta^2+102\theta+15) \\
  &\qquad-2z^2(1129\theta^4+5032\theta^3+7597\theta^2+4773\theta+1083)\\
  &\qquad+2z^3(843\theta^4+2628\theta^3+2353\theta^2+675\theta+6) \\
  &\qquad-z^4(295\theta^4+608\theta^3+478\theta^2+174\theta+26)
     +z^5(\theta + 1)^4
\end{split}
\end{equation*}
This is the Picard-Fuchs differential equation for the complete
intersection of $7$ hyperplanes with the Grassmannian $G(2,7)$ with the
invariants $H^3=42$, $c_2\!\cdot\!H=84$, and $c_3=-98$. (See
\cite{BCKS}.) The singularities are $0$, $3$, $\infty$, and
the three roots $z_1=0.01621\ldots$, $z_2=-0.2139\ldots$, and
$z_3=289.197\ldots$ of $z^3-289z^2-57z+1$. The monodromy around $z=3$
is the identity. The others have the matrix representations
$$
  T_0=\begin{pmatrix}1 & 1 & 0 & 0 \\ 0 & 1 & 0 & 0 \\
  42 & 42 & 1 & 0 \\ 0 & -14 & -1 & 1 \end{pmatrix}, \qquad
  T_{z_1}=\begin{pmatrix} 1 & 0 & 0 & 0 \\ 0 & 1 & 0 & 1 \\
  0 & 0 & 1 & 0 \\ 0 & 0 & 0 & 1\end{pmatrix},
$$
$$
  T_{z_2}=\begin{pmatrix} -13 & 7 & 1 & -2 \\ 28 & -13 & -2 & 4 \\
  -196 & 98 & 15 & -28 \\ 98 & -49 & -7 & 15 \end{pmatrix}, \qquad
  T_{z_3}=\begin{pmatrix} 1 & 0 & 0 & 0 \\ 42 & 1 & 0 & 9 \\
  -196 & 0 & 1 & -42 \\ 0 & 0 & 0 & 1\end{pmatrix}.
$$
Thus, the monodromy group is contained in the subgroup $\Gamma(14,7)$.
\end{section}

\begin{section}{A general approach}
Let
$$
  y^{(n)}+r_{n-1}y^{(n-1)}+\cdots+r_1y^\prime+r_0y=0, \qquad
  r_i\in\C(z),
$$
be a linear differential equation with regular singularities. Then
the monodromy around a singular point $z_0$ with respect to the local
Frobenius basis at $z_0$ is actually very easy to describe, as we
shall see in the following discussion.

Consider the simplest cases where the indicial equation at $z_0$ has
$n$ distinct roots $\lambda_1,\ldots,\lambda_n$ such that
$\lambda_i-\lambda_j\not\in\Z$ for all $i\neq j$. In this case, the
Frobenius basis consists of
$$
  y_j(z)=(z-z_0)^{\lambda_j}f_j(z), \qquad j=1,\ldots,n,
$$
where $f_j(z)$ are holomorphic near $z_0$ and have non-vanishing
constant terms. It is easy to see that the matrix of the monodromy
around $z_0$ with respect to $\{y_j\}$ is simply
$$
  \begin{pmatrix} e^{2\pi i\lambda_1} & 0 & \cdots & 0 \\
   0 & e^{2\pi i\lambda_2} & \cdots & 0 \\
   \vdots & \vdots &  & \vdots \\
   0 & 0 & \cdots & e^{2\pi i\lambda_n}\end{pmatrix}.
$$

Now assume that the indicial equation at $z_0$ has
$\lambda_1,\ldots,\lambda_k$, with multiplicities $e_1,\ldots,e_k$, as
solutions, where $e_1+\cdots+e_k=n$ and $\lambda_i-\lambda_j\not\in\Z$
for all $i\neq j$. Then for each $\lambda_j$, there are $e_j$ linearly
independent solutions
\begin{align*}
  y_{j,0}&=(z-z_0)^{\lambda_j}f_{j,0}, \\
  y_{j,1}&=y_{j,0}\log(z-z_0)+(z-z_0)^{\lambda_j}f_{j,1}, \\
  y_{j,2}&=\frac12 y_{j,0}\log^2(z-z_0)
    +(z-z_0)^{\lambda_j}f_{j,1}\log(z-z_j)+(z-z_0)^{\lambda_j}f_{j,2}, \\
  \vdots & \qquad\qquad\vdots \\
  y_{j,e_j-1}&=(z-z_0)^{\lambda_j}\sum_{h=0}^{e_j-1}\frac1{h!}
    f_{j,e_j-1-h}\log^h(z-z_0),
\end{align*}
where $f_{j,h}$ are holomorphic near $z=z_0$ and satisfy
$f_{j,0}(z_0)=1$ and $f_{j,h}(z_0)=0$ for $h>0$. Since $f_{j,h}$ are
all holomorphic near $z_0$, the analytic continuation along a small
closed curve circling $z_0$ does not change $f_{j,h}$. For other
factors, circling $z_0$ once in the counterclockwise direction results
in
$$
  (z-z_0)^{\lambda_j}\longmapsto e^{2\pi i\lambda_j}(z-z_0)^{\lambda_j}
$$
and
$$
  \log(z-z_0)\longmapsto\log(z-z_0)+2\pi i.
$$
Thus, the behaviors of $y_{j,h}$ under the monodromy around $z_0$ are
governed by
$$
  \begin{pmatrix} y_{j,0} \\ y_{j,1} \\ \vdots \\
    y_{j,e_j-1}\end{pmatrix} \longmapsto
  \begin{pmatrix} \omega_j & 0 & \cdots & 0 \\
   2\pi i\omega_j & \omega_j &   \cdots & 0 \\
   \vdots & \vdots &  & \vdots \\
   \frac{(2\pi i)^{e_j-1}}{(e_j-1)!}\omega_j & \frac{(2\pi
   i)^{e_j-2}}{(e_j-2)!}\omega_j & \cdots & \omega_j\end{pmatrix}
  \begin{pmatrix} y_{j,0} \\ y_{j,1} \\ \vdots \\
    y_{j,e_j-1}\end{pmatrix},
$$
where $\omega_j=e^{2\pi i\lambda_j}$.

When the indicial equation of $z_0$ has distinct roots $\lambda_i$ and
$\lambda_j$ such that $\lambda_i-\lambda_j\in\Z$, there are many
possibilities for the monodromy matrix relative to the Frobenius
basis, but in any case, the matrix still consists of blocks of entries
that take the same form as above.

From the above discussion we see that monodromy matrices with respect
to the local Frobenius bases are very easy to describe. Therefore, to
find monodromy matrices uniformly with respect to a given fixed basis,
it suffices to find the matrix of basis change between the fixed basis
and the Frobenius basis at each singularity. When the differential
equation is hypergeometric, this can be done using the (refined)
standard analytic method, in which we first express the Frobenius
basis at $z=0$ as integrals of Barnes-Mellin type and then use contour
integration to obtain the analytic continuation to a neighborhood of
$z=\infty$. This gives us the monodromy matrices around $z=0$ and
$z=\infty$. Since the monodromy group is generated by these two
matrices, the group is determined.

When the differential equation is not hypergeometric, we are unable to
determine the matrices of basis change precisely. To obtain the
matrices numerically we use the following idea. Let $z_1$ and $z_2$ be
two singularities and $\{y_i\}$ and $\{\tilde y_j\}$, $i,j=1,\ldots,n$,
be their Frobenius bases. Observe that if $y_i=a_{i1}\tilde
y_1+\cdots+a_{in}\tilde y_n$, then we have
$$
  \begin{pmatrix} y_1 & y_1^\prime & \dotsb & y_1^{(n-1)} \\
  y_2 & y_2^\prime & \dotsb & y_2^{(n-1)} \\
  \vdots & \vdots & & \vdots \\
  y_n & y_n^\prime & \dotsb & y_n^{(n-1)} \end{pmatrix}
 =\begin{pmatrix} a_{11} & a_{12} & \dotsb & a_{1n} \\
  a_{21} & a_{22} & \dotsb & a_{2n} \\
  \vdots & \vdots & & \vdots \\
  a_{n1} & a_{n2} & \dotsb & a_{nn} \end{pmatrix}
  \begin{pmatrix} \tilde y_1 & \tilde y_1^\prime & \dotsb
    & \tilde y_1^{(n-1)} \\
  \tilde y_2 & \tilde y_2^\prime & \dotsb & \tilde y_2^{(n-1)} \\
  \vdots & \vdots & & \vdots \\
  \tilde y_n & \tilde y_n^\prime & \dotsb & \tilde y_n^{(n-1)} \end{pmatrix}.
$$
Thus, to determine the matrix $(a_{ij})$ it suffices to evaluate
$y_i^{(k)}$ and $\tilde y_i^{(k)}$ at a common point. To do it
numerically, we expand the Frobenius bases into power series and
assume that the domains of convergence for the power series
have a common point $z_0$. We then truncate and evaluate the series at
$z_0$. This gives us approximation of the matrices of basis changes. We will
discuss some practical issues of this method in Section 5.
\end{section}

\begin{section}{The hypergeometric cases} Throughout this
  section, we fix the branch cut of $\log z$ to be $(-\infty,0]$ so
  that the argument of a complex variable $z$ is between $-\pi$ and
  $\pi$.

  Recall that a hypergeometric function
  $_pF_{p-1}(\alpha_1,\ldots,\alpha_p;\beta_1,\ldots,\beta_{p-1};z)$
  is defined for $\beta_i\neq0,-1,-2,\ldots$ by
  $$
    _pF_{p-1}(\alpha_1,\ldots,\alpha_p;\beta_1,\ldots,\beta_{p-1};z)
   =\sum_{n=0}^\infty\frac{(\alpha_1)_n\ldots(\alpha_p)_n}
    {(1)_n(\beta_1)_n\ldots(\beta_{p-1})_n}z^n,
  $$
  where
  $$
    (\alpha)_n=\begin{cases}\alpha(\alpha+1)\ldots(\alpha+n-1),
    &\text{if }n>0, \\
    1, &\text{if }n=0.
    \end{cases}
  $$
  It satisfies the differential equation
  \begin{equation} \label{temp: hypergeometric}
    [\theta(\theta+\beta_1-1)\ldots(\theta+\beta_{p-1}-1)
   -z(\theta+\alpha_1)\ldots(\theta+\alpha_p)]f=0.
  \end{equation}
  Moreover, it has an integral representation
  $$
    \frac1{2\pi i}\frac{\Gamma(\beta_1)\ldots\Gamma(\beta_{p-1})}
    {\Gamma(\alpha_1)\ldots\Gamma(\alpha_p)}\int_{\mathcal C}
    \frac{\Gamma(s+\alpha_1)\ldots\Gamma(s+\alpha_p)}
    {\Gamma(s+\beta_1)\ldots\Gamma(s+\beta_{p-1})}\Gamma(-s)
    (-z)^s\,ds
  $$
  for $|\arg(-z)|<\pi$, where $\mathcal C$ is any path from $-i\infty$
  to $i\infty$ such that the poles of $\Gamma(-s)$ lie on the right
  of $\mathcal C$ and the poles of $\Gamma(s+a_k)$ lie on the left of
  $\mathcal C$. (See \cite[Chapter 5]{Rainville}.) Then one can obtain the
  analytic continuation of $_pF_{p-1}$ by moving the path of
  integration to the far left of the complex plane and counting the
  residues arising from the process. It turns out that this method can
  be generalized.

  \begin{Lemma} Let $m$ be the number of $1$'s among $\beta_k$. Set
  $$
    F(h,z)=\sum_{n=0}^\infty\frac{(\alpha_1+h)_n\ldots(\alpha_p+h)_n}
    {(1+h)_n(\beta_1+h)_n\ldots(\beta_{p-1}+h)_n}z^{n+h}.
  $$
  Then, for $j=0,\ldots,m$, the functions
  $$
    \frac{\partial^j}{\partial h^j}F(h,z)\Big|_{h=0}
  $$
  are solutions of (\ref{temp: hypergeometric}). Moreover, if $|\arg
  (-z)|<\pi$ and $h$ is a small quantity such that $\alpha_k+h$ are not
  zero or negative integers, then $F(h,z)$ has the integral
  representation
  \begin{equation*}
  \begin{split}
    F(h,z)&=-\frac{z^h}{2\pi i}\frac{\Gamma(\beta_1+h)\ldots
     \Gamma(\beta_{p-1}+h)\Gamma(1+h)}
    {\Gamma(\alpha_1+h)\ldots\Gamma(\alpha_p+h)} \\
   &\qquad\int_{\mathcal C}
    \frac{\Gamma(s+\alpha_1+h)\ldots\Gamma(s+\alpha_p+h)}
    {\Gamma(s+\beta_1+h)\ldots\Gamma(s+\beta_{p-1}+h)\Gamma(s+1+h)}
    \frac\pi{\sin\pi s}(-z)^s\,ds,
  \end{split}
  \end{equation*}
  where $\mathcal C$ is any path from $-i\infty$ to $i\infty$ such
  that the integers $0,1,2,\ldots$ lies on the right of $\mathcal C$
  and the poles of $\Gamma(s+a_k+h)$ lie on the left of $\mathcal C$.
  \end{Lemma}

  \begin{proof} The first part of the lemma is just a specialization
  of the Frobenius method (see \cite{Ince}) to the hypergeometric cases.
  We have
  \begin{equation*}
  \begin{split}
   &\theta(\theta+\beta_1-1)\ldots(\theta+\beta_{p-1}-1)F(h,z)
     =h(h+\beta_1-1)\ldots(h+\beta_{p-1}-1)z^h \\
   &\qquad+\sum_{n=1}^\infty\frac{(\alpha_1+h)_n\ldots(\alpha_p+h)}
    {(1+h)_{n-1}(\beta_1+h)_{n-1}\ldots(\beta_{p-1}+h)_{n-1}}z^{n+h}
  \end{split}
  \end{equation*}
  and
  \begin{equation*}
  \begin{split}
    &z(\theta+\alpha_1)\ldots(\theta+\alpha_p)F(h,z)
    =\sum_{n=0}^\infty\frac{(\alpha_1+h)_{n+1}\ldots(\alpha_p+h)_{n+1}}
    {(1+h)_n(\beta_1+h)_n\ldots(\beta_{p-1}+h)_n}z^{n+1+h}
  \end{split}
  \end{equation*}
  It follows that
  \begin{equation*}
  \begin{split}
   &[\theta(\theta+\beta_1-1)\ldots(\theta+\beta_{p-1}-1)
   -z(\theta+\alpha_1)\ldots(\theta+\alpha_p)]F(h,z) \\
   &\qquad=h(h+\beta_1-1)\ldots(h+\beta_{p-1}-1)z^h.
  \end{split}
  \end{equation*}
  If the number of $1$'s among $\beta_k$ is $m$, then the first
  non-vanishing term of the Taylor expansion in $h$ of the last
  expression is $h^{m+1}$. Consequently, we see that
  $$
    \frac{\partial^j}{\partial h^j}F(h,z)\Big|_{h=0}
  $$
  are solutions of (\ref{temp: hypergeometric}) for $j=0,\ldots,m$.

  The proof of the second part about the integral representation is
  standard. We refer the reader to Chapter 5 of \cite{Rainville}.
  \end{proof}

  We now prove Theorem 1. Here we will only discuss the cases
$$(A,B)=(1/2,1/2),\, (1/3,1/3),\, (1/4,1/2),\,\,\mbox{and}\,\,(1/6,1/3),$$
  representing the four classes whose indicial equations at $z=\infty$
  have one root with multiplicity $4$, two distinct roots, each of
  which has multiplicity $2$, one repeated root and two other distinct
  roots, and four distinct roots, respectively. The other cases can be
  proved in the same fashion.

  \begin{proof}[Proof of the case $(A,B)=(1/6,1/3)$] Let $h$ denote a
  small real number, and let $F(h,z)$ be defined as in Lemma 1 with
  $p=4$, $\alpha_1=1/6$, $\alpha_2=1/3$, $\alpha_3=2/3$,
  $\alpha_4=5/6$, and $\beta_k=1$ for all $k$. Then, by Lemma 1, the
  functions
  $$
    y_j(z)=\frac1{j!}\frac{\partial^j}{\partial h^j}(C^{-h}F(h,Cz)),
    \qquad j=0,\ldots,3,
  $$
  are solutions of
  $$
    \theta^4-11664z(\theta+1/6)(\theta+1/3)(\theta+2/3)(\theta+5/6),
  $$
  where $C=11664$. In fact, by considering the contribution of the
  first term, we see that these four functions make up the Frobenius
  basis at $z=0$.

  We now express $C^{-h}F(h,Cz)$ using Lemma 1. By the Gauss
  multiplication theorem we have
  \begin{equation*}
  \begin{split}
   &\Gamma(s+1/6)\Gamma(s+1/3)\Gamma(s+2/3)\Gamma(s+5/6) \\
   &\qquad=\frac{\prod_{k=1}^6\Gamma(s+k/6)}{\Gamma(s+1/2)\Gamma(s+1)}
   =\frac{(2\pi)^{5/2}6^{-1/2-6s}\Gamma(6s+1)}
    {(2\pi)^{1/2}2^{-1/2-2s}\Gamma(2s+1)}.
  \end{split}
  \end{equation*}
  Thus, restricting $z$ to the lower half-plane $-\pi<\arg z<0$, by
  Lemma 1, we may write
  \begin{equation*}
  \begin{split}
    C^{-h}F(h,Cz)&=-\frac{z^h}{2\pi i}\frac{\Gamma(1+h)^4\Gamma(1+2h)}
    {\Gamma(1+6h)} \\
   &\qquad\times\int_{\mathcal C}\frac{\Gamma(6s+1+6h)}
    {\Gamma(s+1+h)^4\Gamma(2s+1+2h)}\frac\pi{\sin\pi s}
    e^{\pi is}z^s\,ds,
  \end{split}
  \end{equation*}
  where $\mathcal C$ is the vertical line $\Re s=-1/12$.
  Now move the line of integration to $\Re s=-13/12$. This is
  justified by the fact that the integrand tends to $0$ as $\Im s$
  tends to infinity. The integrand has four simple poles $s=-n/6-h$,
  $n=1,2,4,5$, between these two lines. The residues are
  $$
    \frac{(-1)^{n-1}}{6\Gamma(n)}\frac{\pi e^{-\pi i(n/6+h)}}
    {\Gamma(1-n/6)^4\Gamma(1-n/3)\sin\pi(n/6+h)}z^{-n/6-h}.
  $$
  Thus, we see that the analytic continuation of $C^{-h}F(h,z)$ to
  $|z|>1$ with $-\pi<\arg z<0$ is given by
  \begin{equation*}
  \begin{split}
    C^{-h}F(h,z)=\sum_{n=1,2,4,5}a_n B_n(h)z^{-n/6}+
   (\text{higher order terms in }1/z),
  \end{split}
  \end{equation*}
  where
  $$
    a_n=\frac{(-1)^n\pi e^{-\pi in/6}}{6\Gamma(n)\Gamma(1-n/6)^4
    \Gamma(1-n/3)}, \qquad
    B_n(h)=\frac{\Gamma(1+h)^4\Gamma(1+2h)e^{-\pi ih}}
    {\Gamma(1+6h)\sin\pi(n/6+h)}.
  $$
  On the other hand, since the local exponents at
  $z=\infty$ are $1/6$, $1/3$, $2/3$, and $5/6$, the Frobenius basis
  at $z=\infty$ consists of
  $$
    \tilde y_n(z)=z^{-n/6}g_n(1/z), \qquad n=1,2,4,5,
  $$
  where $g_n=1+\cdots$ are functions holomorphic at $0$. It follows
  that for $z$ with $-\pi<\arg z<0$
  $$
    y_j(z)=\frac1{j!}\sum_{n=1,2,4,5}a_n B_n^{(j)}(h)\tilde y_n(z).
  $$
  Set $f_j(z)=y_j(z)/(2\pi i)^j$ for $j=0,\ldots,3$ and $\tilde
  f_n=a_n\tilde y_n/\sin(n\pi /6)$ for $n=1,2,4,5$. Then using the
  evaluation
  $$
    \Gamma^\prime(1)=-\gamma, \qquad
    \Gamma^{\prime\prime}(1)=\gamma^2+\zeta(2), \qquad
    \Gamma^{\prime\prime\prime}(1)=-\gamma^3-3\zeta(2)\gamma-2\zeta(3),
  $$
  we find
  $$
    \begin{pmatrix} f_3 \\ f_2 \\ f_1 \\ f_0 \end{pmatrix}
  =M\begin{pmatrix} \tilde f_1 \\ \tilde f_2 \\ \tilde f_4 \\
    \tilde f_5\end{pmatrix},
  $$
  where
  $$
   M=\begin{pmatrix}
    \eta-i\omega/4 &\eta+5\sqrt 3i\omega^2/36
   &\eta+5\sqrt 3i\omega^4/36 &\eta-i\omega^5/4 \\
    -5/12-i\omega/2 & 1/4-i\omega^2/2\sqrt 3&
    1/4-i\omega^4/2\sqrt 3 & -5/12-i\omega^5/2 \\
    i\omega & i\omega^2/\sqrt 3 & i\omega^4/\sqrt 3 & i\omega^5 \\
    1 & 1 & 1 & 1 \end{pmatrix}
  $$
  with
  $$
    \omega=e^{\pi i/6}, \qquad \eta=\frac{68\zeta(3)}{(2\pi i)^3}.
  $$
  Now let $P$ be the path traveling along the real axis with $\arg
  z=-\pi+$ from $z=-2$ to $-\infty$ and then coming back along the
  real axis with $\arg z=\pi-$ to $z=-2$. The monodromy effect on
  $\tilde y_n(z)=z^{-n/6}g_n(1/z)$ is
  $$
    \tilde y_n(z)\longmapsto\tilde y_n(e^{2\pi i}z)
    =e^{-2\pi in/6}\tilde y_n(z).
  $$
  Therefore, the matrix representation of the monodromy along $P$
  relative to the ordered basis $\{f_3,f_2,f_1,f_0\}$ is
  $$
    T_\infty=M\begin{pmatrix} \omega^{-2} & 0 & 0 & 0 \\
     0 & \omega^{-4} & 0 & 0 \\ 0 & 0 & \omega^4 & 0 \\
     0 & 0 & 0 & \omega^2\end{pmatrix} M^{-1}.
  $$
  Now the path $P$ is equivalent to that of circling once around
  $z=1/C$ and then once around $z=0$, both in the counterclockwise
  direction. Therefore, if we denote by $T_0$ and $T_{1/C}$ the
  monodromy matrices relative the basis $\{f_3,f_2,f_1,f_0\}$ around
  $z=0$ and $z=1/C$, respectively, then we have
  $$
    T_\infty=T_{1/C}T_0.
  $$
  Since $T_0$ is easily seen to be
  $$
    T_0=\begin{pmatrix} 1 & 1 & 1/2 & 1/6 \\
     0 & 1 & 1 & 1/2 \\ 0 & 0 & 1 & 1 \\ 0 & 0 & 0 & 1 \end{pmatrix},
  $$
  we find
  $$
    T_{1/C}=M\begin{pmatrix} \omega^{-2} & 0 & 0 & 0 \\
     0 & \omega^{-4} & 0 & 0 \\ 0 & 0 & \omega^4 & 0 \\
     0 & 0 & 0 & \omega^2\end{pmatrix} M^{-1}T_0^{-1}
   =\begin{pmatrix} 1+a & 0 & ab/d & a^2/d \\
    -b & 1 & -b^2/d & -ab/d \\ 0 & 0 & 1 & 0 \\
    -d & 0 & -b & 1-a \end{pmatrix},
  $$
  where
  $$
    a=-\frac{204}{(2\pi i)^3}\zeta(3), \qquad b=\frac 74,
    \qquad d=3.
  $$
  Comparing these numbers with the invariants, we find the matrix
  $T_{1/C}$ indeed takes the form (\ref{raw matrices}) specified in
  the statement of Theorem 1. This proves the case $(A,B)=(1/6,1/3)$.
  \end{proof}

  \begin{proof}[Proof of the case $(A,B)=(1/4,1/2)$] Apply Lemma 1
  with $p=4$, $\alpha_1=1/4$, $\alpha_2=3/4$, $\alpha_3=\alpha_4=1/2$,
  $\beta_k=1$ for all $k$, and set $C=1024$. Then
  $$
    y_j(z)=\frac1{j!}\frac{\partial^j}{\partial h^j}(C^{-h}F(h,Cz)),
    \qquad j=0,\ldots,3,
  $$
  form the Frobenius basis for
  $$
    \theta^4-1024z(\theta+1/4)(\theta+3/4)(\theta+1/2)^2.
  $$
  Assuming that $-\pi<\arg z<0$, we have
  \begin{equation*}
  \begin{split}
    C^{-h}F(h,Cz)&=-\frac{z^h}{2\pi i}\frac{\Gamma(1+h)^6}
    {\Gamma(1+2h)\Gamma(1+4h)} \\
   &\qquad\times\int_{\mathcal C}\frac{\Gamma(4s+1+4h)\Gamma(2s+1+2h)}
    {\Gamma(s+1+h)^6}\frac\pi{\sin\pi s}
    e^{\pi is}z^s\,ds,
  \end{split}
  \end{equation*}
  where $\mathcal C$ is the vertical line $\Re s=-1/8$. The integrand
  has simple poles at $-k-h-1/4$ and $-k-h-3/4$, and double poles at
  $-k-h-1/2$ for $k=0,1,2,\ldots$. The residues at $s=-h-n/4$,
  $n=1,3$, are $a_nC_n(h)z^{-h-n/4}$, where
  $$
   a_n=(-1)^{(n+1)/2}\frac{\pi\Gamma(1/2)e^{-\pi in/4}}{4\Gamma(1-n/4)^6},
   \qquad C_n(h)=\frac{e^{-\pi ih}}{\sin\pi(h+n/4)}
  $$
  At $s=-h-1/2$ we have
  \begin{equation*}
  \begin{split}
    &\frac{\Gamma(4s+1+4h)\Gamma(2s+1+2h)}{\Gamma(s+1+h)^6} \\
    &\qquad=-\frac1{8\Gamma(1/2)^6}(s+h+1/2)^{-2}-\frac{3\log
      2+1}{2\Gamma(1/2)^6}(s+h+1/2)^{-1}+\cdots,
  \end{split}
  \end{equation*}
  $$
    \frac\pi{\sin\pi s}=-\frac\pi{\cos\pi h}
   +\pi^2\frac{\sin\pi h}{\cos^2\pi h}(s+h+1/2)+\cdots,
  $$
  and
  $$
    e^{\pi is}z^s=z^{-1/2-h}e^{-\pi i(h+1/2)}(1+(\pi i+\log z)(s+h+1/2)
   +\cdots).
  $$
  Thus, the residue at $s=-h-1/2$ is
  $$
    \frac{\pi e^{-\pi i(h+1/2)}}{8\Gamma(1/2)^6\cos\pi h}\left(
    \pi i+\log z+12\log 2+4-\pi\frac{\sin\pi h}{\cos\pi h}\right)
    z^{-h-1/2}.
  $$
  Set
  $$
    a_2=-\frac{\pi e^{-\pi i/2}}{8\Gamma(1/2)^6}, \quad
    C_2(h)=\frac{e^{-\pi ih}}{\cos\pi h}, \quad
    C_2^\ast(h)=C_2(h)(\pi i+12\log 2+4-\pi\tan\pi h).
  $$
  We find
  \begin{equation*}
  \begin{split}
    C^{-h}F(h,Cz)&=-a_1B_1(h)z^{-1/4}-a_2B_2(h)z^{-1/2}\log z
    -a_2B^\ast_2(h)z^{-1/2} \\
    &\qquad\qquad-a_3B_3(h)z^{-3/4}+(\text{higher order terms in }1/z),
  \end{split}
  \end{equation*}
  where
  $$
    B_n(h)=\frac{\Gamma(1+h)^6}{\Gamma(1+2h)\Gamma(1+4h)}C_n(h),
    \qquad B^\ast_2(h)
   =\frac{\Gamma(1+h)^6}{\Gamma(1+2h)\Gamma(1+4h)}C^\ast_2(h).
  $$
  Let $y_j(z)$, $j=0,\ldots,3$, be the Frobenius basis at $z=0$, and
  \begin{equation*}
  \begin{split}
    \tilde y_1(z)=z^{-1/4}(1+\cdots), & \qquad
    \tilde y_3(z)=z^{-3/4}(1+\cdots), \\
    \tilde y_2^\ast(z)=z^{-1/2}(1+\cdots), & \qquad
    \tilde y_2(z)=(\log z+g(1/z))\tilde y_2^\ast(z)
  \end{split}
  \end{equation*}
  be the Frobenius basis at $\infty$, where $g(t)$ is a function
  holomorphic and vanishing at $t=0$. Set $f_j(z)=y_j(z)/(2\pi i)^j$
  for $j=0,\ldots,3$, $\tilde f_n(z)=-a_n\tilde y_n(z)/\sin\pi(n/4)$
  for $n=1,2,3$, and $\tilde f_2^\ast(z)=-a_2\tilde y^\ast_2(z)$.
  Using the fact that
  $$
    y_j(z)=\frac1{j!}\frac{\partial^j}{\partial h^j}(C^{-h}F(h,Cz)),
  $$
  we find
  $$
    \begin{pmatrix} f_3 \\ f_2 \\ f_1 \\ f_0 \end{pmatrix}
   =\begin{pmatrix}
    \eta+(1-i)/48& \eta-5/48 & -5\mu/12+\pi i\eta+4\mu\eta
      & \eta+(1+i)/48\\ 
    (1-6i)/24 & 7/24 & \pi i/24+7\mu/6 & (1+6i)/24 \\
    (i-1)/2 & -1/2 & -2\mu & -(i+1)/2 \\
    1 & 1 & \pi i+4\mu & 1 \end{pmatrix}
    \begin{pmatrix} \tilde f_1 \\ \tilde f_2 \\ \tilde f_2^\ast \\
    \tilde f_3\end{pmatrix}
  $$
  where
  $$
    \mu=3\log 2+1, \qquad \eta=\frac{22\zeta(3)}{(2\pi i)^3}.
  $$
  Let $P$ be the path from $z=-1$ with argument $-\pi$ to $-\infty$
  and then back to $z=-1$ with argument $\pi$. The monodromy matrix
  for $P$ relative to the ordered basis $\{\tilde f_1,\tilde f_2,
  \tilde f_2^\ast, \tilde f_3\}$ is
  $$
    \begin{pmatrix} -i & 0 & 0 & 0 \\ 0 & -1 & -2\pi i & 0 \\
    0 & 0 & -1 & 0 \\ 0 & 0 & 0 & i\end{pmatrix}.
  $$
  Thus, the matrix with respect to the ordered basis
  $\{f_3,f_2,f_1,f_0\}$ is
  $$
    T_\infty=\begin{pmatrix}
    1-8\eta & 1-8\eta & 1/2-19\eta/3 & 1/6-11\eta/3+8\eta^2 \\
    -7/3 & -4/3 & -61/72 & -41/72+7\eta/3 \\
    0 & 0 & 1 & 1 \\ -8 & -8 & -19/3 & -8/3+8\eta\end{pmatrix}.
  $$
  Finally, it is easy to see that the monodromy around $z=0$ with
  respective to $\{f_3,f_2,f_1,f_0\}$ is
  $$
    T_0=\begin{pmatrix} 1 & 1 & 1/2 & 1/6 \\ 0 & 1 & 1 & 1/2 \\
    0 & 0 & 1 & 1 \\ 0 & 0 & 0 & 1\end{pmatrix},
  $$
  and after a short computation we find that the monodromy
  $T_{1/C}=T_\infty T_0^{-1}$ around $z=1/C$ indeed takes the form
  claimed in the statement of Theorem 1.
  \end{proof}

  \begin{proof}[Proof of the case $(A,B)=(1/3,1/3)$] Let $z$ be a
  complex number with $-\pi<\arg z<0$. By the same argument as before.
  We find that the Frobenius basis $\{y_j\}$ at $z=0$ can be expressed
  as
  $$
    y_j(z)=\frac1{j!}\frac{\partial^j}{\partial h^j}(C^{-h}F(h,Cz)),
  $$
  where $C=729$ and
  \begin{equation*}
  \begin{split}
    C^{-h}F(h,Cz)&=-\frac{z^h}{2\pi i}\frac{\Gamma(1+h)^6}
    {\Gamma(1+3h)^2}\int_{\mathcal C}\frac{\Gamma(3s+1+3h)^2}
    {\Gamma(s+1+h)^6}\frac\pi{\sin\pi s}
    e^{\pi is}z^s\,ds.
  \end{split}
  \end{equation*}
  Here $h$ is assumed to be a real number and $\mathcal C$ denotes the
  vertical line $\Re s=-1/6$. Set
  $$
    a_n=-\frac{\pi e^{-\pi in/3}}{9\Gamma(1-n/3)^6}, \qquad n=1,\ 2.
  $$
  The residues at $z=-1/3-h$ and $z=-2/3-h$ are
  $$
    a_1(\pi i+\log z+9\log 3-\pi\sqrt 3+\cot\pi(1/3+h))
    z^{-1/3-h}e^{-\pi ih}
  $$
  and
  $$
    a_2(\pi i+\log z+9\log 3+\pi\sqrt 3+6+\cot\pi(2/3+h))
    z^{-2/3-h}e^{-\pi ih},
  $$
  respectively. Let
  $$
    B_n(h)=\frac{\Gamma(1+h)^6}{\Gamma(1+3h)^2}
    \frac{e^{-\pi ih}}{\sin\pi(n/3+h)}, \qquad n=1,\ 2,
  $$
  and
  $$
    B_1^\ast(h)=B_1(h)(\pi i+9\log 3-\pi\sqrt 3+\pi\cot\pi(1/3+h)),
  $$
  $$
    B_2^\ast(h)=B_2(h)(\pi i+9\log 3+\pi\sqrt 3+6+\pi\cot\pi(2/3+h)).
  $$
  Then we have
  $$
    C^{-h}F(h,Cz)=-\sum_{n=1}^2a_n(B_n(h)z^{-n/3}\log z+
    B_n^\ast(h)z^{-n/3})+\text{(higher order terms)}.
  $$
  Now the Taylor expansions of $B_n(h)$ and $B_n^\ast(n)$ are
  $$
    \sin\frac\pi 3B_1\left(\frac h{2\pi i}\right)
   =1+\frac{i\omega}{\sqrt3}h
     -\left(\frac{i\omega}{2\sqrt 3}+\frac1{12}\right)h^2
     +\left(\frac{i\omega}{12\sqrt 3}+\eta\right)h^3+\cdots,
  $$
  $$
    \sin\frac\pi 3B_2\left(\frac h{2\pi i}\right)
   =1+\frac{i\omega^2}{\sqrt3}h
     -\left(\frac{i\omega^2}{2\sqrt 3}+\frac1{12}\right)h^2
     +\left(\frac{i\omega^2}{12\sqrt 3}+\eta\right)h^3+\cdots,
  $$
  \begin{equation*}
  \begin{split}
    \sin\frac\pi 3B_1^\ast\left(\frac h{2\pi i}\right)
  &=\left(\mu_1+\frac{2\pi\omega}{\sqrt 3}\right)
   +\left(\frac{i\mu_1\omega}{\sqrt 3}+\frac{2\pi i\omega}3\right)h
   -\left(\frac{i\mu_1\omega}{2\sqrt 3}+\frac{\mu_1}{12}
     +\frac{\pi\omega}{2\sqrt 3}\right)h^2 \\
  &\qquad+\left(\mu_1\eta+\frac{i\mu_1\omega}{12\sqrt 3}
     +\frac{2\pi\eta\omega}{\sqrt 3}-\frac{\pi i\omega}6\right)h^3
   +\cdots,
  \end{split}
  \end{equation*}
  \begin{equation*}
  \begin{split}
    \sin\frac\pi 3B_2^\ast\left(\frac h{2\pi i}\right)
  &=\left(\mu_2+\frac{2\pi\omega^2}{\sqrt 3}+6\right)
   +\left(\frac{i\mu_2\omega^2}{\sqrt 3}-\frac{2\pi i\omega^2}3
     +2i\omega^2\sqrt 3\right)h \\
  &\qquad
   -\left(\frac{i\mu_2\omega^2}{2\sqrt 3}+\frac{\mu_2}{12}
     +\frac{\pi\omega^2}{2\sqrt 3}-\omega^2-\frac32\right)h^2 \\
  &\qquad+\left(\mu_2\eta+\frac{i\mu_2\omega^2}{12\sqrt 3}
     +\frac{2\pi\eta\omega^2}{\sqrt 3}+\frac{\pi i\omega^2}6
     +6\eta+\frac{i\omega^2}{2\sqrt 3}\right)h^3
   +\cdots,
  \end{split}
  \end{equation*}
  where
  $$
    \omega=e^{\pi i/3}, \qquad \eta=\frac{16\zeta(3)}{(2\pi i)^3},
    \qquad\mu_1=9\log 3-\pi\sqrt 3, \qquad\mu_2=9\log 3+\pi\sqrt 3.
  $$
  From these we can deduce the matrix of basis change between the
  Frobenius basis
  $$
    f_j(z)=\frac1{(2\pi i)^jj!}\frac{\partial^j}{\partial h^j}
      C^{-h}F(h,Cz), \qquad j=0,\ldots,3
  $$
  and the basis
  \begin{equation*}
  \begin{split}
    \tilde f_1^\ast(z)=-a_1z^{-1/3}(1+\cdots), &\qquad
    \tilde f_2^\ast(z)=-a_2z^{-2/3}(1+\cdots), \\
    \tilde f_1=(\log z+g_1(1/z))\tilde f_1^\ast(z), &\qquad
    \tilde f_2=(\log z+g_2(1/z))\tilde f_2^\ast(z),
  \end{split}
  \end{equation*}
  where $g_1(t)$ and $g_2(t)$ are functions holomorphic and vanishing
  at $t=0$. The monodromy matrix around $\infty$ with respect to the
  ordered basis $\{\tilde f_1,\tilde f_1^\ast, \tilde f_2, \tilde
  f_2^\ast\}$ is easily seen to be
  $$
    \begin{pmatrix}
    e^{-2\pi i/3} & 2\pi ie^{-2\pi i/3} & 0 & 0 \\
    0 & e^{-2\pi i/3} & 0 & 0 \\
    0 & 0 & e^{2\pi i/3} & 2\pi ie^{2\pi i/3} \\
    0 & 0 & 0 & e^{2\pi i/3}
    \end{pmatrix}.
  $$
  By the same argument as before, we find that the monodromy matrix
  with respect to the basis $\{f_3,f_2,f_1,f_0\}$ indeed takes the
  form claimed in the statement. This proves the case
  $(A,B)=(1/3,1/3)$.
  \end{proof}

  \begin{proof}[Proof of the case $(A,B)=(1/2,1/2)$] Let $z$ be a
  complex number such that $-\pi<\arg z<0$. We find that the Frobenius
  basis $\{y_j\}$ at $z=0$ can be expressed as
  $$
    y_j(z)=\frac1{j!}\frac{\partial^j}{\partial h^j}(C^{-h}F(h,Cz)),
  $$
  where $C=256$ and
  \begin{equation*}
  \begin{split}
    C^{-h}F(h,Cz)&=-\frac{z^h}{2\pi i}\frac{\Gamma(1+h)^8}
    {\Gamma(1+2h)^4}\int_{\mathcal C}\frac{\Gamma(2s+1+2h)^4}
    {\Gamma(s+1+h)^8}\frac\pi{\sin\pi s}
    e^{\pi is}z^s\,ds.
  \end{split}
  \end{equation*}
  Here $h$ is assumed to be a small real number and $\mathcal C$
  denotes the vertical line $\Re s=-1/4$. The integrand has
  quadruple poles at $s=-k-1/2-h$ for non-positive integers $k$.
  Moving the line of integration to $\Re s=-3/4$ and computing the
  residue at $s=-1/2-h$, we see that
  $$
    C^{-h}F(h,Cz)=a_1\sum_{n=0}^3\frac{B_n(h)}{n!}
    z^{-1/2}(\log z)^n+(\text{higher order terms in }1/z),
  $$
  where
  $$
    a_1=\frac{\pi e^{-\pi i/2}}{16\Gamma(1/2)^8}, \quad
    B_3(h)=\frac{\Gamma(1+h)^8e^{-\pi ih}}{\Gamma(1+2h)^4\cos\pi h},
    \quad B_2(h)=B_3(h)(\mu-\pi\tan\pi h),
  $$
  $$
    B_1(h)=B_3(h)\left(-\frac 76\pi^2+\frac{\mu^2}2-\pi\mu\tan\pi h
     +\pi^2\sec^2\pi h\right),
  $$
  and
  \begin{equation*}
  \begin{split}
    B_0(h)&=B_3(h)\Big(\frac{\mu^3}6-\frac\pi 2\mu^2\tan\pi h
      +(\sec^2\pi h-7/6)\pi^2\mu  \\
    &\qquad\qquad\qquad+(5/6-\sec^2\pi h)\pi^3\tan\pi h
      +8\zeta(3)\Big),
  \end{split}
  \end{equation*}
  where $\mu=16\log2+\pi i$. Let
  $$
    \tilde f_0(z)=z^{-1/2}(1+\cdots), \qquad
    \tilde f_1(z)=\frac1{2\pi i}(\log z+g_1(1/z))\tilde f_0(z),
  $$
  $$
    \tilde f_2(z)=\frac1{(2\pi i)^2}(\log^2z/2+g_1(1/z)\log
    z+g_2(1/z))\tilde f_0(z),
  $$
  $$
    \tilde f_3(z)=\frac1{(2\pi i)^3}
   (\log^3z/6+g_1(1/z)\log^2z/2+g_2(1/z)\log z+g_3(1/z))\tilde f_0(z)
  $$
  be the Frobenius basis at $z=\infty$ with $g_n(0)=0$. Using the
  evaluation
  $$
    \Gamma^\prime(1)=-\gamma, \qquad
    \Gamma^{\prime\prime}(1)=\frac{\pi^2}{12}+\frac{\gamma^2}2, \qquad
    \Gamma^{\prime\prime\prime}(1)=-\frac13\zeta(3)
   -\frac{\pi^2\gamma}{12}-\frac{\gamma^3}6,
  $$
  we can find the analytic continuation of the Frobenius at $z=0$ in
  terms of $\tilde f_n(z)$. Now the monodromy around $\infty$ relative
  to the basis $\{\tilde f_3(z),\tilde f_2(z),\tilde f_1(z), \tilde
  f_0(z)\}$ is
  $$
    \begin{pmatrix} -1 & -1 & -1/2 & -1/6 \\
     0 & -1 & -1 & -1/2 \\ 0 & 0 & -1 & -1 \\ 0 & 0 & 0 & -1
    \end{pmatrix}.
  $$
  From this we can determine the monodromy matrix around $z=1/C$ with
  respect to the Frobenius basis at $z=0$. We find that the result
  agrees with the general pattern depicted in Theorem 1, although the
  detailed computation is too complicated to be presented here.
  \end{proof}

  Of course, there is no reason why our approach should be applicable
  only to order $4$ cases. Consider the hypergeometric differential
  equations of the form
  \begin{equation} \label{order 5}
    L:\theta^5-z(\theta+1/2)(\theta+A)(\theta+1-A)(\theta+B)(\theta+1-B).
  \end{equation}
  The cases $(A,B)=(1/2,1/2)$, $(1/4,1/2)$,
  $(1/6,1/4)$, $(1/4,1/3)$, $(1/6,1/3)$, and
  $(1/8,3/8)$ have been used by Guillera \cite{Guillera1,Guillera2}
  to construct series representations for $1/\pi^2$. Applying the above
  method, we determine the monodromy of these differential equations
  in the following theorem whose proof will be omitted.

  \begin{Theorem} Let $L$ be one of the differential equations in
  (\ref{order 5}). Let $y_i$, $i=0,\ldots,4$, be the Frobenius basis
  at $0$. Then the monodromy matrices around $z=0$ and $z=1/C$ with
  respect to the ordered basis $\{y_4/(2\pi i)^4,y_3/(2\pi i)^3,
  y_2/(2\pi i)^2, y_1/(2\pi i), y_0\}$ are
  $$
    \MF{11{1/2}{1/6}{1/24}}{011{1/2}{1/6}}{0011{1/2}}
    {00011}{00001}, \quad
    \MF{{a^2}0{-ab}{(1-a^2)x}{-b^2/2}}
    {{-c^2x/2}1{-acx}{c^2x^2/2}{-(1-a^2)x}}
    {{-ac}0{1-2a^2}{acx}{-ab}}{00010}
    {{-c^2/2}0{-ac}{c^2x/2}{a^2}},
  $$
  respectively, where $x$ is an integer multiple of $\zeta(3)/(2\pi
  i)^3$, $a$ and $c$ are positive real numbers such that $a^2$, $ac$,
  and $c^2$ are rational numbers, and $b$ is a real number satisfying
  $a^2+bc=1$. The exact values of $a$, $c$, and $x^\prime=(2\pi
  i)^3x/\zeta(3)$ are given in the following table.
  $$ \arraycolsep3pt
  \begin{array}{||cc|ccc||} \hline
  A & B & a^2 & c^2 & x^\prime \\ \hline
  1/2 & 1/2 & 25/36 & 64 & 10\\
  1/2 & 1/4 & 8/9 & 32 & 24 \\
  1/4 & 1/6 & 289/288 & 8 & 80 \\
  1/3 & 1/4 & 27/32 & 24 & 28 \\
  1/3 & 1/6 & 75/64 & 12 & 70 \\
  1/8 & 3/8 & 529/288 & 8 & 150 \\
  \hline
  \end{array}
  $$
  \end{Theorem}
\end{section}

\begin{section}{Differential equations of Calabi-Yau type} The
  Picard-Fuchs differential equations for families of Calabi-Yau
  threefolds known in literature have the common features that
  \begin{enumerate}
  \item[(a)] the singular points are all regular,
  \item[(b)] the indicial equation at $z=0$ has $0$ as its only solution,
  \item[(c)] the indicial equation at one of the singularities has
  solutions $0,1,1,2$, corresponding to a conifold singularity,
  \item[(d)] the unique holomorphic solution $y$ around $0$ with $y(0)=1$
    has integral coefficients in its power series expansion,
  \item[(e)] the solutions $\lambda_1\le\lambda_2\le\lambda_3\le\lambda_4$
    of the indicial equation at $t=\infty$ are positive rational
    numbers and satisfy $\lambda_1+\lambda_4=\lambda_2+\lambda_3=r$
    for some $r\in\Q$, and the characteristic polynomial of the
    monodromy around $t=\infty$ is a product of cyclotomic
    polynomials.
  \item[(f)] the coefficients $r_i(z)$ of the differential equation satisfy
   $$
     r_1=\frac12r_2r_3-\frac18r_3^3+r_2^\prime-\frac34r_3^\prime r_3
    -\frac12 r_3^{\prime\prime},
   $$
  \item[(g)] the instanton numbers are integers.
  \end{enumerate}
  In \cite{AESZ} a fourth order linear differential equation satisfying
  all conditions except (c) is said to be of {\it Calabi-Yau type}.
  Using various techniques, Almkvist and etc. found more than
  $300$ such equations. (See Section 5 of \cite{Almkvist-Zudilin}
  for an overview of strategies of finding Calabi-Yau equations. The
  paper also contains a ``superseeker'' that tabulates the known
  Calabi-Yau equations, sorted according to the instanton numbers.)
  Among them, there are $178$ equations that have singularities with
  exponents $0,1,1,2$. It is speculated that all such equations should
  have geometric origin.

  In \cite{Enckevort-Straten} van Enckevort and van Straten
  numerically determined the monodromy for these $178$ equations.
  They were able to find rational bases for $145$ of them, among which
  there are $64$ cases that are integral. Their method goes as
  follows. Let $z_1,\ldots,z_k$ be the singularities of a Calabi-Yau
  differential equation. They first chose a reference point $p$ and
  piecewise linear loops each of which starts from $p$ and encircles
  exactly one of $z_i$. Then the problem of determining analytic
  continuation becomes that of solving several initial value problems
  in sequences. This was done numerically using the {\tt dsolve}
  function in {\tt Maple}. Then they used the crucial observation that
  the Jordan form for the monodromy around a conifold singularity is
  unipotent of index one to find a rational basis. Finally, assuming
  that (\ref{DM representation}) and (\ref{conifold period}) hold for
  general differential equations of Calabi-Yau type, conjectural
  values of geometric invariants can be read off.

  Here we present a different method of computing monodromy based
  on the approach described in Section 3. Let $0=z_0,z_1,\ldots,z_n$ be
  the singular points of a Calabi-Yau differential equation, and
  assume that $f_{i,k}$, $i=0,\ldots,n$, $k=1,\ldots,4$ form the
  Frobenius bases at $z_i$. According to Section 3, to find the matrix
  of basis change between $\{f_{i,k}\}$ and $\{f_{j,k}\}$, we only
  need to evaluate $f_{i,k}^{(m)}$ and $f_{j,k}^{(m)}$ at a common
  point $\zeta$ where the power series expansions of the functions
  involved all converge. In practice, the choice of $\zeta$ is
  important in order to achieve required precision in a reasonable
  amount of time.

  Let $R_i$ denote the radius of convergence of the power series
  expansions of the Frobenius basis at $z_i$. In general, $R_i$ is
  equal to the distance from $z_i$ to the nearest singularity $z_j\neq
  z_i$, meaning that if we truncate the power series expansion of
  $f_{i,k}$ at the $n$th term, the resulting error is
  $$
    O_\epsilon\left((1+\epsilon)^n\frac{|\zeta-z_i|^n}{R_i^n}\right).
  $$
  Of course, the $O$-constants depend on the differential equation
  and $z_i$. Since we do not have any control over them, in practice
  we just choose $\zeta$ in a way such that
  $$
    \frac{|\zeta-z_i|}{R_i}=\frac{|\zeta-z_j|}{R_j}.
  $$
  If this does not yield needed precision, we simply replace $n$ by a
  larger integer and do the computation again.
  \medskip

  \noindent{\bf Example.} Consider
  $$
    \theta^4-5(5\theta+1)(5\theta+2)(5\theta+3)(5\theta+4).
  $$
  The singularities are $z_0=0$, $z_1=1/3125$, and $z_2=\infty$.
  The radii of convergence for the Frobenius bases at $0$ and $1/3125$
  are both $1/3125$. Thus, to find the matrix of basis change, we
  expand the Frobenius bases, say, for $30$ terms, and evaluate the
  Frobenius bases and their derivatives at $\zeta=1/6250$. Then we use
  the idea in Section 3 to compute the monodromy matrix around $z_1$
  with respective to the Frobenius basis at $0$. We find that the
  computation agrees with (\ref{raw matrices}) in Theorem 1 up to $7$
  digits.

  \medskip

  The above method works quite well if the singularities of a
  differential equation are reasonably well spaced. However, it occurs
  quite often that a Calabi-Yau differential equation has a cluster of
  singular points near $0$, and a couple of singular points that are
  far away. For example, consider Equation \#19
  \begin{equation*}
  \begin{split}
    &529\theta^4-23z(921\theta^4 + 2046\theta^3 + 1644\theta^2 +
      621\theta + 92) \\
    &\qquad - z^2(380851\theta^4 + 1328584\theta^3 +
      1772673\theta^2 + 1033528\theta + 221168) \\
    &\qquad - 2z^3(475861\theta^4 + 1310172\theta^3 + 1028791\theta^2
      + 208932\theta - 27232) \\
    &\qquad - 68z^4(8873\theta^4 + 14020\theta^3 + 5139\theta^2 -
      1664\theta - 976) \\
    &\qquad + 6936z^5(\theta + 1)^2(3\theta + 2)(3\theta + 4).
  \end{split}
  \end{equation*}
  The singularities are $z_0=0$, $z_1=1/54$,
  $z_2=(11-5\sqrt 5)/2=-0.090\ldots$, $z_3=-23/34$, and
  $z_4=(11+5\sqrt 5)/2=11.09\ldots$. In order to determine the
  monodromy matrix around $z_4$, we need to compute the matrix of
  basis change between the Frobenius basis at $1/54$ and that at
  $z_4$. The radius of convergence for the Frobenius basis at $1/54$
  is $1/54$, while that at $z_4$ is $z_4-1/54=11.07\ldots$. Even if we
  choose $\zeta$ optimally, we still need to expand the Frobenius
  bases for thousands of terms in order to achieve a precision of a few
  digits. In such situations, we can choose several points lying
  between the two singularities, compute bases for each of them, and
  then use the same idea as before to determine the matrices of basis
  change.

  Take Equation 19 above as an example. We choose $w_k=(1+3^k)/54$ and
  $\zeta_k=(1+3^k/2)/54$ for $k=0,\ldots,5$. The radius of convergence
  for the basis at $w_k$ is $3^k/54$. Thus, evaluating the first $n$
  terms of the power series expansions at $\zeta_k$ and $\zeta_{k+1}$
  will result in an error of
  $$
    O_\epsilon((1/2+\epsilon)^n),
  $$
  which is good enough in practice.

  Using the above ideas we computed the monodromy groups of the
  differential equations of Calabi-Yau type that have at least one
  conifold singularity.\footnote{\bf We have written a Maple program
  for the computation. It is available upon request.} Our result shows
  that if a differential equation comes from geometry, then the
  monodromy matrix around one of the conifold singularities with
  respect to the Frobenius basis at the origin takes the form
  (\ref{raw matrices}). We then conjugate the monodromy matrices by
  the matrix (\ref{CY basis}) and find that the other matrices are
  also in $\Sp(4,\Q)$. We now tabulate the results for the equations
  coming from geometry in the following table. Note that the notations
  $\Gamma(d_1,d_2)$ and $\Gamma(d_1,d_2,d_3)$, $d_2,d_3|d_1$,
  represent the congruence subgroups
  \begin{equation*}
  \begin{split}
    \Gamma(d_1,d_2)&=\left\{\gamma\in\Sp(4,\Z):\gamma\equiv
    \begin{pmatrix} 1 & \ast & \ast & \ast \\
     0 & \ast & \ast & \ast \\ 0 & 0 & 1 & 0 \\
     0 & \ast & \ast & \ast \end{pmatrix} \mod d_1\right\} \\
     &\qquad\qquad\bigcap\left\{\gamma\in\Sp(4,\Z):\gamma\equiv
    \begin{pmatrix} 1 & \ast & \ast & \ast \\
     0 & 1 & \ast & \ast \\ 0 & 0 & 1 & 0 \\
     0 & 0 & \ast & 1 \end{pmatrix} \mod d_2\right\}
  \end{split}
  \end{equation*}
  and
  \begin{equation*}
  \begin{split}
    \Gamma(d_1,d_2,d_3)
  &=\Big\{(a_{ij})\in\Sp(4,\Q): a_{ij}\in\Z\ \forall(i,j)\neq(1,3),\
    a_{13}\in\frac1{d_3}\Z, \\
  &\qquad\qquad a_{21}, a_{31}, a_{41}, a_{32}, a_{34}\equiv 0\mod d_1, \\
  &\qquad\qquad a_{42}\equiv0, \quad a_{22}, a_{44}\equiv 1\mod d_2, \\
  &\qquad\qquad a_{11},a_{33}\equiv 1\mod\frac{d_1}{d_3}\Big\}.
  \end{split}
  \end{equation*}
  Note also that since the matrix
  $$
    \begin{pmatrix} 1 & 0 & 0 & 0\\ 0 & 1 & 0 & 1\\ 0 & 0 & 1 & 0\\
    0 & 0 & 0 & 1 \end{pmatrix}
  $$
  is always in the monodromy groups, it is not listed in the table.
  The reader should be mindful of this omission.
  $$
    \extrarowheight3pt
    \arraycolsep4pt
    \begin{array}{|c||c|c|c||l|c||c|} \hline
    \# & H^3 & c_2\!\cdot\! H\!\! & c_3
       & \text{Generators} & \text{in} & \text{Ref} \\ \hline
    15 & 18 & 72 & -162
       & \M{1100}{0100}{{18}{18}10}{0{-9}{-1}1}
         \M{{-2}3{1/2}{-1}}{6{-5}{-1}2}{{-18}{18}4{-6}}
           {{18}{-18}{-3}7}
       & \Gamma(6,3,2) & \text{\cite{Batyrev-Straten}} \\ \hline
    16 & 48 & 96 & -128
       & \M{1100}{0100}{{48}{48}10}{0{-16}{-1}1}
         \M{{-5}2{1/4}{-1}}{{24}{-7}{-1}4}{{-144}{48}7{-24}}
           {{48}{-16}{-2}9}
        & \Gamma(24,8,4) & \text{\cite{Batyrev-Straten}} \\ \hline
    17 & 90 & 108 & -90
       & \begin{array}{l}
         \!\!\M{1100}{0100}{{90}{90}10}{0{-24}{-1}1}
         \M{{-17}3{1/3}1}{{-54}{10}13}
           {{-972}{162}{19}{54}}{{162}{-27}{-3}{-8}} \\
         \!\!\M{{-11}3{1/3}{-1}}{{36}{-8}{-1}3}
           {{-432}{108}{13}{-36}}{{108}{-27}{-3}{10}}
         \end{array}
       & \Gamma(18,6,3) & \text{\cite{Batyrev-Straten}} \\ \hline
    18 & 40 & 88 & -128
       & \M{1100}{0100}{{40}{40}10}{0{-14}{-1}1}
         \M{{-5}4{1/2}{-1}}{{12}{-7}{-1}2}{{-72}{48}7{-12}}
           {{48}{-32}{-4}9}
        & \Gamma(4,2,2) & \text{\cite{Batyrev-Straten}} \\ \hline
    19 & 46 & 88 & -106
       & \begin{array}{l}
         \!\!\M{1100}{0100}{{46}{46}10}{0{-15}{-1}1}
         \M{{-6}4{1/2}{-1}}{{14}{-7}{-1}2}
           {{-98}{56}8{-14}}{{56}{-32}{-4}9} \\
         \!\!\M{{-45}{12}2{-6}}{{138}{-35}{-6}{18}}
           {{-1058}{276}{47}{-138}}{{276}{-72}{-12}{37}}
         \end{array}
       &\Gamma(2,2,2) & \text{\cite{Batyrev-Straten}} \\ \hline
    20 & 54 & 72 & -18
       & \M{1100}{0100}{{54}{54}10}{0{-15}{-1}1}
         \M{7{-1}{-1/6}1}{{-6}10{-2}}{{126}{-18}{-2}{24}}{{-36}61{-5}}
       & \Gamma(6,3,6) & \text{\cite{Batyrev-Straten}} \\ \hline
    21 & 80 & 104 & -88
       & \begin{array}{l}
         \!\!\M{1100}{0100}{{80}{80}10}{0{-22}{-1}1}
         \M{{-11}5{1/2}{-1}}{{24}{-9}{-1}2}{{-288}{120}{13}{-24}}
           {{120}{-50}{-5}{11}} \\
         \!\!\M{{-19}4{1/2}{-2}}{{80}{-15}{-2}8}{{-800}{160}{21}{-80}}
           {{160}{-32}{-4}{17}}
         \end{array}
       & \Gamma(8,2,2) & \text{\cite{Batyrev-Straten}} \\ \hline
    22 & 70 & 100 & -100
       & \begin{array}{l}
         \!\!\M{1100}{0100}{{70}{70}10}{0{-20}{-1}1}
         \M{1000}{{10}102}{{-50}01{-10}}{0001} \\
         \!\!\M{{-9}5{1/2}{-1}}{{20}{-9}{-1}2}
         {{-200}{100}{11}{-20}}{{100}{-50}{-5}{11}}
         \end{array}
       & \Gamma(10,10,2) & \text{\cite{Batyrev-Straten}} \\ \hline
    23 & 96 & 96 & -32
       & \M{1100}{0100}{{96}{96}10}{0{-24}{-1}1}
         \M{9{-1}{-1/8}1}{{-8}10{-2}}
         {{288}{-32}{-3}{40}}{{-64}81{-7}}
       & \Gamma(8,8,8) & \text{\cite{Batyrev-Straten}} \\ \hline
    \end{array}
  $$
  $$
    \extrarowheight3pt
    \arraycolsep4pt
    \begin{array}{|c||c|c|c||l|c||c|} \hline
    24 & 15 & 66 & -150
       & \M{1100}{0100}{{15}{15}10}{0{-8}{-1}1}
         \M{{-5}51{-2}}{{12}{-9}{-2}4}{{-36}{30}7{-12}}
         {{30}{-25}{-5}{11}}
       & \Gamma(3,1) & \text{\cite{BCKS}} \\ \hline
    25 & 20 & 68 & -120
       & \M{1100}{0100}{{20}{20}10}{0{-9}{-1}1}
         \M{{-7}51{-2}}{{16}{-9}{-2}4}{{-64}{40}9{-16}}
         {{40}{-25}{-5}{11}}
       & \Gamma(4,1) & \text{\cite{BCKS}} \\ \hline
    26 & 28 & 76 & -116
       & \M{1100}{0100}{{28}{28}10}{0{-11}{-1}1}
         \M{{-9}61{-2}}{{20}{-11}{-2}4}{{-100}{60}{11}{-20}}
         {{60}{-36}{-6}{13}}
       & \Gamma(4,1) & \text{\cite{BCKS}} \\ \hline
    27 & 42 & 84 & 98
       & \begin{array}{l}
         \!\!\M{1100}{0100}{{42}{42}10}{0{-14}{-1}1}
         \M{1000}{{42}109}{{-196}01{-42}}{0001} \\
         \!\!\M{{-13}71{-2}}{{28}{-13}{-2}4}{{-196}{98}{15}{-28}}
         {{98}{-49}{-7}{15}}
         \end{array}
       & \Gamma(14,7) & \text{\cite{BCKS,Rodland}} \\ \hline
    28 & 42 & 84 & -96
       & \M{1100}{0100}{{42}{42}10}{0{-14}{-1}1}
         \M{{-41}{12}2{-6}}{{126}{-35}{-6}{18}}
         {{-882}{252}{43}{-126}}{{252}{-72}{-12}{37}}
       & \Gamma(42,2) & \text{\cite{BCKS}} \\ \hline
   186 & 57 & 90 & -84
       & \begin{array}{l}
       \!\!\M{1100}{0100}{{57}{57}10}{0{-17}{-1}1}
         \M{{-53}{12}2{-6}}{{162}{-35}{-6}{18}}{{-1458}{324}{55}{-162}}
           {{324}{-72}{-12}{37}} \\
       \!\!\M{{-17}81{-2}}{{36}{-15}{-2}4}{{-324}{144}{19}{-36}}
           {{144}{-64}{-8}{17}}
         \end{array}
       & \Gamma(3,1) & \text{\cite{Tjotta}} \\ \hline
    \end{array}
  $$

  In the second table we list a few equations whose monodromy matrices
  with respect to our bases have integers as entries. Note that
  the numbers $H^3$, $c_2\!\cdot\! H$, and $c_3$ are all conjectural,
  obtained from evaluation of the monodromy around a singularity of
  conifold type. Note, again, that the matrix
  $$
    \begin{pmatrix} 1 & 0 & 0 & 0\\ 0 & 1 & 0 & 1\\ 0 & 0 & 1 & 0\\
    0 & 0 & 0 & 1 \end{pmatrix}
  $$
  is omitted from the table.
  $$
    \extrarowheight3pt
    \arraycolsep4pt
    \begin{array}{|c||c|c|c||l|c|} \hline
    \# & H^3 & c_2\!\cdot\! H\!\! & c_3
       & \text{Generators} & \text{in} \\ \hline
    29 & 24 & 72 & -116
       & \M{1100}{0100}{{24}{24}10}{0{-10}{-1}1}
         \M{{-47}{20}4{-10}}{{120}{-49}{-10}{25}}
         {{-576}{240}{49}{-120}}{{240}{-100}{-20}{51}}
       & \Gamma(24,2) \\ \hline
    33 & 6 & 36 & -72
       & \M{1100}{0100}{6610}{0{-4}{-1}1}
         \M{1000}{2102}{{-2}01{-2}}{0001}
       & \Gamma(2,2) \\ \hline
    42 & 32 & 80 & -116
       & \M{1100}{0100}{{32}{32}10}{0{-12}{-1}1}
         \M{{-15}61{-3}}{{48}{-17}{-3}9}{{-256}{96}{17}{-48}}
           {{96}{-36}{-6}{19}}
       & \Gamma(16,4) \\ \hline
    \end{array}
  $$
  $$
    \extrarowheight3pt
    \arraycolsep4pt
    \begin{array}{|c||c|c|c||l|c|} \hline
    51 & 10 & 64 & -200
       & \M{1100}{0100}{{10}{10}10}{0{-7}{-1}1}
         \M{{-3}51{-2}}{8{-9}{-2}4}{{-16}{20}5{-8}}{{20}{-25}{-5}{11}}
       & \Gamma(2,1) \\ \hline
    63 & 5 & 62 & {-310}
       & \M{1100}{0100}{5510}{0{-6}{-1}1}
         \M{{-1}51{-2}}{4{-9}{-2}4}{{-4}{10}3{-4}}{{10}{-25}{-5}{11}}
       & \Gamma(1,1) \\ \hline
    73 & 9 & 30 & 12
       & \M{1100}{0100}{9910}{0{-4}{-1}1}
         \M{2001}{3{-2}{-1}0}{0323}{{-3}00{-2}}
       & \Gamma(3,1) \\ \hline
    99 & 13 & 58 & -120
       & \M{1100}{0100}{{13}{13}10}{0{-7}{-1}1}
         \M{{-5}41{-2}}{{12}{-7}{-2}4}{{-36}{24}7{-12}}
          {{24}{-16}{-4}9}
       & \Gamma(1,1) \\ \hline
   100 & 36 & 72 & -72
       & \begin{array}{l}
         \!\!\M{1100}{0100}{{36}{36}10}{0{-12}{-1}1}
         \M{1000}{{12}104}{{-36}01{-12}}{0001} \\
         \!\!\M{{-11}61{-2}}{{24}{-11}{-2}4}{{-144}{72}{13}{-24}}
           {{72}{-36}{-6}{13}}
         \end{array}
       & \Gamma(12,12) \\ \hline
   101 & 25 & 70 & -100
       & \begin{array}{l}
         \!\!\M{1100}{0100}{{25}{25}10}{0{-10}{-1}1}
         \M{{-19}{10}2{-4}}{{40}{-19}{-4}8}{{-200}{100}{21}{-40}}
           {{100}{-50}{-10}{21}} \\
         \!\!\M{1000}{{60}10{16}}{{-225}01{-60}}{0001}
         \end{array}
       & \Gamma(5,5) \\ \hline
   109 & 7 & 46 & -120
       & \M{1100}{0100}{7710}{0{-5}{-1}1}
         \M{{-3}31{-2}}{8{-5}{-2}4}{{-16}{12}5{-8}}{{12}{-9}{-3}7}
       & \Gamma(1,1) \\ \hline
   117 & 12 & 36 & -32
       & \begin{array}{l}
         \!\!\M{1100}{0100}{{12}{12}10}{0{-5}{-1}1}
         \M{1000}{4104}{{-4}01{-4}}{0001} \\
         \!\!\M{{-59}{21}9{18}}{{-120}{43}{18}{36}}{{-400}{140}{61}{120}}
           {{140}{-49}{-21}{-41}} \\
         \end{array}
       & \Gamma(4,1) \\ \hline
   118 & 10 & 40 & -50
       & \begin{array}{l}
         \!\!\M{1100}{0100}{{10}{10}10}{0{-5}{-1}1}
         \M{1000}{{30}10{18}}{{-50}01{-30}}{0001} \\
         \!\!\M{{-19}{10}4{-8}}{{40}{-19}{-8}{16}}{{-100}{50}{21}{-40}}
           {{50}{-25}{-10}{21}}
         \end{array}
       & \Gamma(10,5) \\ \hline
   185 & 36 & 84 & -120
       & \M{1100}{0100}{{36}{36}10}{0{-13}{-1}1}
         \M{{-11}71{-2}}{{24}{-13}{-2}4}{{-144}{84}{13}{-24}}
           {{84}{-49}{-7}{15}}
       & \Gamma(12,1) \\ \hline
    \end{array}
  $$
\end{section}
\vfill\newpage

\centerline{\bf Acknowledgments}
\medskip

The second author (Yifan Yang) would like to thank Wadim Zudilin for
drawing his attention to the monodromy problems and for many
interesting and fruitful discussions. This whole research project
started out as the second author and Zudilin's attempt to give a
rigorous and uniform proof of Guillera's $1/\pi^2$ formulas
\cite{Guillera1,Guillera2} in a way analogous to the modular-function
approach in \cite{Chan-Chan-Liu}. (See also \cite{Yang}.) For this
purpose, it was natural to consider the monodromy of the fifth-order
hypergeometric differential equations, and hence it led the second
author to consider monodromy of general differential equations of
Calabi-Yau type. The second author would also like to thank Duco van
Straten for his interest in this project and for clarifying some
questions about the differential equations.

During the preparation of this paper (at the final stage), the third
author (N. Yui) was a visiting researcher at Max-Planck-Institut f\"ur
Mathematik Bonn in May and June 2006. Her visit was supported by
Max-Planck-Institut. She thanks Don Zagier, Andrey Todorov and Wadim
Zudilin for their interest in this project and discussions and
suggestions on the topic discussed in this paper.  She is especially
indebted to Cord Erdenberger of University of Hannover for his
supplying the index calculation for the congruence subgroup
$\Gamma(d_1,d_2)$ in $\Sp(4,\Z)$.

The final version was prepared at IHES in January 2007 where the third
author was a visiting member.  For the preparation of the fnial version,
the comments and suggestions of the referee were very helpful
as well as discussions with Maxim Kontsevich. We thank them for their help.

\vfill\newpage
\centerline{\bf Appendix -- The index of $\Gamma(d_1,d_2)$ in $\Sp(4,\Z)$}
\medskip

\centerline{by {\sc Cord Erdenberger}}
\medskip
\centerline{
  Institut f\"ur Mathematik, Universit\"at Hannover}
\centerline{
  Welfengarten 1, 30167 Hannover, Germany}
\centerline{\tt
  E-mail: erdenber@math.uni-hannover.de}
\medskip

In this appendix we will show that the groups $\Gamma(d_1,d_2)$ are
indeed congruence subgroups in $\Sp(4,\Z)$ and provide a formula for
their index.
\medskip

Recall that for $n\in\N$ the principal congruence subgroup of level
$n$ is defined by
$$
\Gamma(n):=\left\{\,M\in \Sp(4,\Z)\,|\, N\equiv \mathbf{I}_4\pmod n\,\right\}.
$$
It is the kernel of the map from $\Sp(4,\Z)$ to $\Sp(4,\Z/n\Z)$ given
by reduction modulo $n$ and thus a normal subgroup in $\Sp(4,\Z)$. It
is a well--known fact that this map is surjective and hence the
sequence
$$
  \I_4\to \Gamma(n)\hookrightarrow \Sp(4,\Z)\to \Sp(4,\Z/n\Z)\to \I_4
$$
is exact. So the index of $\Gamma(n)$ in $\Sp(4,\Z)$ is just the order
of $\Sp(4,\Z/n\Z)$ which is known to be
$$
[\Sp(4,\Z):\Gamma(n)]=|\Sp(4,\Z/n\Z)|=n^{10}\prod (1-p^{-2})(1-p^{-4}),
$$
where the product runs over all primes $p$ such that $p|n$. 
\medskip

For $d_1,d_2 \in \N$, define
$$
\widetilde{\Gamma}_1(d_1):=\left\{\,M\in \Sp(4,\Z)\,:\,M\equiv\begin{pmatrix}
1 & * & * & * \\
0 & * & * & * \\
0 & 0 & 1 & 0 \\
0 & * & * & *
\end{pmatrix}\, \mod{d_1}\,\right\},
$$
$$
\widetilde{\Gamma}_2(d_2):=\left\{\, M\in \Sp(4,\Z)\,:\,
M\equiv\begin{pmatrix}
1 & * & * & * \\
0 & 1 & * & * \\
0 & 0 & 1 & 0 \\
0 & 0 & * & 1 \end{pmatrix}\,\mod{d_2}\,\right\}
$$
and set
$$
  \Gamma(d_1,d_2):=\widetilde{\Gamma}_1(d_1)\cap\widetilde{\Gamma}_2(d_2).
$$
Note that 
$$
  \Gamma(d_1)\subset\widetilde{\Gamma}_1(d_1)\,\,\mbox{and}\,\,
  \Gamma(d_2)\subset\widetilde{\Gamma}_2(d_2).
$$
Hence
$$
  \Gamma(d)=\Gamma(d_1)\cap\Gamma(d_2)\subset\widetilde{\Gamma}_1(d_1)
  \cap\widetilde{\Gamma}_2(d_2)=\Gamma(d_1,d_2),
$$
where $d$ is the least common multiple of $d_1$ and $d_2$. This shows
that $\Gamma(d_1,d_2)$ is a congruence subgroup, i.e. it contains a
principal congruence subgroup as a normal subgroup of finite
index. Moreover, this implies that $\Gamma(d_1,d_2)$ has finite index
in $\Sp(4,\Z)$ and an upper bound is given by the index of $\Gamma(d)$
as given above.
\medskip

We will from now on restrict to the case relevant to this paper,
namely $d_2|d_1$. Then $\Gamma(d_1,d_2)$ is in fact a subgroup of
$\widetilde{\Gamma}_1(d_1)$, namely
$$
\Gamma(d_1,d_2)=\left\{\begin{pmatrix} a_{11} & a_{12} & a_{13} & a_{14} \\
  a_{21} & a_{22} & a_{23} & a_{24} \\
  a_{31} & a_{32} & a_{33} & a_{34} \\
  a_{41} & a_{42} & a_{43} & a_{44} \end{pmatrix} \in
\widetilde{\Gamma}_1(d_1) :
  \begin{pmatrix}
        a_{22} & a_{24} \\
        a_{42} & a_{44}
  \end{pmatrix}
  \equiv\begin{pmatrix}
  1 & * \\
  0 & 1 
\end{pmatrix} \bmod{d_2} \right\} .
$$

To obtain a formula for the index of this group in $\Sp(4,\Z)$, we
first calculate the index of $\widetilde{\Gamma}_1(d_1)$.
Note that
$$
  \widetilde{\Gamma}_1(d_1)/\Gamma(d_1) < \Sp(4,\Z)/\Gamma(d_1)
  \simeq \Sp(4,\Z/d_1\Z)
$$
and hence
$$
 [\Sp(4,\Z):\widetilde{\Gamma}_1(d_1)]=[\Sp(4,\Z/d_1\Z):
  \widetilde{\Gamma}_1(d_1)/\Gamma(d_1)].
$$
The quotient $\widetilde{\Gamma}_1(d_1)/\Gamma(d_1)$ considered as a
subgroup of $\Sp(4,\Z/d_1\Z)$ via the above isomorphism is given by 
$$\widetilde{\Gamma}_1(d_1)/\Gamma(d_1) \simeq \left\{\, M\in
  \Sp(4,\Z/d_1\Z)\,:\,
M=\begin{pmatrix}
1 & * & * & * \\
0 & * & * & * \\
0 & 0 & 1 & 0 \\
0 & * & * & *
\end{pmatrix}\,\right\} .
$$
An element of this group has the following form
$$M=\begin{pmatrix} 
1 & a_{12} & a_{13} & a_{14} \\
0 & \alpha & a_{23} & \beta \\
0 & 0 & 1 & 0 \\
0 & \gamma & a_{43} & \delta \end{pmatrix}.$$
Let $\mathbb{J}_4:=\begin{pmatrix} 0 & - \mathbf{I}_2 \\
                                   \mathbf{I_2} & 0 \end{pmatrix}$.
The symplectic relation that $^tM\,\mathbb{J}_4\,M=\mathbb{J}_4$ then
implies that $\begin{pmatrix} \alpha & \beta \\ 
  \gamma & \delta \end{pmatrix}\in
\mathrm{SL}_2(\Z/d_1\Z)$. Furthermore it gives rise to the following
linear system:
$$ 
\begin{matrix} a_{12} &+ & \alpha \, a_{43}&-& \gamma \, a_{23}&=& 0\\
               a_{14}&+&\beta \, a_{43}&-& \delta \, a_{23}&=& 0
             \end{matrix}$$
Writing this in matrix form, we have
$$\begin{pmatrix} \alpha & \gamma \\
                  \beta & \delta \end{pmatrix} \begin{pmatrix} -a_{43}
                  \\ a_{23} \end{pmatrix} 
=\begin{pmatrix} a_{12} \\ a_{14} \end{pmatrix}.$$
If we choose $a_{12}, a_{13}, a_{14}$ freely, the above linear system
                  has a unique solution
$a_{23}, a_{43}$ as $\begin{pmatrix} \alpha & \beta \\ \gamma & \delta
                  \end{pmatrix}$ is in 
                  $\mathrm{SL}_2(\Z/d_1\Z)$. This shows that
$$
  |\widetilde{\Gamma}_1(d_1)/\Gamma(d_1)|=d_1^3\cdot
  |\mathrm{SL}_2(\Z/d_1\Z)|=d_1^6\prod(1-p^{-2})
$$
where the product runs over all primes $p$ dividing $d_1$. 
So we have the index formula
$$
  [\Sp(4,\Z):\widetilde{\Gamma}_1(d_1)]=[\Sp(4,\Z/d_1\Z):
  \widetilde{\Gamma}_1(d_1)/\Gamma(d_1)]
=d_1^4 \prod_{p|d_1} (1-p^{-4}).
$$
\medskip

Now we are ready to calculate the index of $\Gamma(d_1,d_2)$ in
$\Sp(4,\Z)$. Since we assume that $d_2|d_1$, we have the following
chain of subgroups:
$$
  \Gamma(d_1) < \Gamma(d_1,d_2) < \widetilde{\Gamma}_1(d_1) < \Sp(4,\Z)
$$

Note that
$$
  [\widetilde{\Gamma}_1(d_1):\Gamma(d_1,d_2)] =
  [\widetilde{\Gamma}_1(d_1)/\Gamma(d_1):\Gamma(d_1,d_2)/\Gamma(d_1)]
$$
and by our above description this is just the index of the group
$$
  \left\{ M \in \mathrm{SL}_2(\Z/d_1\Z)\,:\, M \equiv \begin{pmatrix}
  1 & * \\
  0 & 1 
  \end{pmatrix} \mod{d_2} \right\}
$$
in $\mathrm{SL}_2(\Z/d_1\Z)$. An easy calculation shows that this
index is equal to
$$d_2^2\prod_{p|d_2} (1-p^{-2}).$$

Putting all these together, we get
\begin{align*}
  [\Sp(4,\Z):\Gamma(d_1,d_2)] & =
  [\Sp(4,\Z):\widetilde{\Gamma}_1(d_1)] \cdot
  [\widetilde{\Gamma}_1(d_1):\Gamma(d_1,d_2)] \\
  & = d_1^4 \prod_{p|d_1} (1-p^{-4}) \; d_2^2 \prod_{p|d_2} (1-p^{-2}) .
\end{align*}
\medskip

We summarize the above calculation to obtain
\medskip

\noindent{\bf Theorem.} The group $\Gamma(d_1,d_2)$ is a congruence
subgroup in $\Sp(4,\Z)$ and its index is given by
$$
  |\Sp(4,\Z):\Gamma(d_1,d_2)|=d_1^4\prod_{p|d_1}(1-p^{-4}) \;
  d_2^2\prod_{p|d_2}(1-p^{-2}) .
$$

In fact, we can do a similiar calculation without the assumption that
$d_2|d_1$ and obtain the same formula as given above where one has to
replace $d_1$ with the least common multiple of $d_1$ and $d_2$.

\bibliography{plain}

\end{document}